\documentclass[DIV=14,letterpaper]{scrartcl}

\date{}

\usepackage{tikz}
\usetikzlibrary{calc}
\usetikzlibrary{matrix}
\usepackage{amsmath,amssymb,amsfonts,amsthm,subfig}
\usepackage[pdftex]{hyperref} 

\usepackage{graphicx}

\newcommand{\V}[1]{\mbox{\boldmath $ #1 $}}

\def \p{\partial}

\newcommand{\bey}{\begin{eqnarray}}
\newcommand{\eey}{\end{eqnarray}}

\newcommand{\beq}{\begin{equation}}
\newcommand{\eeq}{\end{equation}}
\newtheorem{thm}{\hspace{6mm}Theorem}[section]

\newtheorem{lem}{\hspace{6mm}Lemma}[section]

\newtheorem{rem}{\hspace{6mm}Remark}[section]
\newcommand{\proofend}{\mbox{ }\hfill \raisebox{.4ex}{\framebox[1ex]{}}}

\title{The cutoff method for the numerical computation of nonnegative solutions of
parabolic PDEs with application to anisotropic diffusion and lubrication-type equations}

\author{Changna Lu\thanks{College of Math and Statistics,
Nanjing University of Information Science and Technology,
Nanjing, Jiangsu 210044, China. ({\tt luchangna@nuist.edu.cn})}
\and Weizhang Huang\thanks{
Department of Mathematics, the University of Kansas, Lawrence, KS 66045,
U.S.A. ({\tt huang@math.ku.edu})}
\and Erik S. Van Vleck\thanks{
Department of Mathematics, the University of Kansas, Lawrence, KS 66045,
U.S.A. ({\tt evanvleck@math.ku.edu})}
}

\begin{document}
\vskip 1cm
\maketitle

\begin{abstract}
The cutoff method,  which cuts off the values of a function less than a given number,
is studied for the numerical computation of nonnegative solutions
of parabolic partial differential equations. A convergence analysis is given for a broad class of
finite difference methods combined with cutoff for linear parabolic equations.
Two applications are investigated, linear anisotropic diffusion problems satisfying the setting
of the convergence analysis and nonlinear lubrication-type equations
for which it is unclear if the convergence analysis applies.
The numerical results are shown to be consistent with the theory
and in good agreement with existing results in the literature.
The convergence analysis and applications demonstrate that the cutoff method is an effective
tool for use in the computation of nonnegative solutions. Cutoff can also be used with other
discretization methods such as collocation, finite volume, finite element, and spectral
methods and for the computation of positive solutions.
\end{abstract}



\noindent{\bf Abbreviated title.}
Cutoff method for computation of nonnegative solutions

\section{Introduction}

Many physical phenomena involve variables such as the density and concentration of a material
that take only nonnegative values. The mathematical reflection of this property is that
the partial differential equations (PDEs) modeling the phenomena admit nonnegative
solutions. For the numerical solution of those PDEs, it is crucial that numerical schemes
preserve the solution nonnegativity and produce physically meaningful numerical solutions.

A closely related concept is the maximum principle. Preserving the maximum principle
is equivalent to preserving the solution nonnegativity for linear problems \cite{TW08} but generally
the former is more difficult than the latter. It is known (e.g., see \cite{CR73}) that a conventional numerical method,
such as a finite difference (FD), finite volume (FV), or finite element (FE) method,
generally does not preserve the maximum principle and can produce negative undershoot in the solution
for diffusion problems, especially those with heterogeneous anisotropic diffusion coefficients.
Considerable effort has been made on developing numerical schemes satisfying the maximum principle; 
e.g., see \cite{BKK08, BE04, Cia70, CR73, KK09, KKK07, KL95, Let92, Sto86, SF73, WaZh11, XZ99}
for steady-state isotropic diffusion problems and
\cite{DDS04, GL09, GYK05, Hua10, KSS09, LH10, LSS07, LSSV07, LS08, LuHuQi2012, MD06, SH07}
for steady-state anisotropic diffusion problems. Particularly, Ciarlet and Raviart \cite{CR73} show that the linear
FE approximation of an elliptic isotropic diffusion problem satisfies the maximum principle
when the mesh is simplicial and nonobtuse. The result is generalized to anisotropic diffusion problems
by Huang and his coworkers in \cite{Hua10,LH10, LuHuQi2012}.
On the other hand, less progress has been made for time dependent problems.
For example, Fujii \cite{Fuj73} shows that the linear FE approximation of the heat equation satisfies the maximum
principle when the mesh is simplicial and acute and the time step size is bounded above
by a bound proportional to the squared maximal element size and below by a bound
proportional to the squared minimal element in-diameter.
He also shows that when the lumped mass matrix is used, the maximum principle holds without requiring the time step
size to be bounded from below.
Fujii's results are extended to more general isotropic diffusion problems by Farag{\'o} et al. \cite{FH2001,FHK05,FHK2011} and to anisotropic diffusion problems by Li and Huang \cite{LiHu2012}.
Farag{\'o} and Horv{\'a}th \cite{FH06} study the relations between the maximum principle, nonnegativity
preservation, and maximum norm contractivity for linear parabolic equations and their finite difference
and Galerkin finite element discretizations.
Thom{\'e}e and Wahlbin \cite{TW08} consider more general parabolic
PDEs and show that
the maximum principle cannot hold for the conventional semidiscrete FE problem
and Fujii's conditions on the mesh are essentially sharp for the lumped mass matrix.
Le Potier \cite{LePot05,LePot09b} (also see Lipnikov et al. \cite{LSSV07})
proposes two nonlinear FV schemes for linear anisotropic diffusion problems and shows that they are
second order in space and satisfy monotonicity \cite{LePot05} or the maximum principle \cite{LePot09b}.

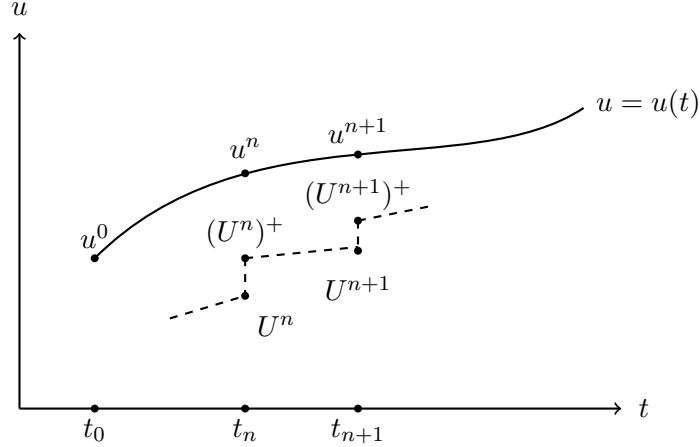
\begin{figure}[tb]
\centering
\begin{tikzpicture}[scale = 1]
\draw [thick, ->] (0,1.5) -- (8,1.5);
\draw [thick, ->] (0,1.5) -- (0,6.5);
\draw[right] (8.1,1.5) node {$t$};
\draw[above] (0,6.6) node {$u$};

\draw[thick] (1,3.5) .. controls (3,5.5) and (6, 4.5) .. (7.5, 5.5);
\draw[right] (7.53, 5.55) node {$u=u(t)$};
\fill (canvas cs:x=1cm,y=3.5cm) circle (1.5pt);
\draw[above] (1,3.5) node {$u^0$};
\fill (canvas cs:x=3cm,y=4.63cm) circle (1.5pt);
\draw[above] (3,4.7) node {$u^n$};
\fill (canvas cs:x=4.5cm,y=4.88cm) circle (1.5pt);
\draw[above] (4.5,4.9) node {$u^{n+1}$};

\fill (canvas cs:x=3cm,y=3cm) circle (1.5pt);
\draw[left] (3.8,2.6) node {${U}^n$};

\draw[thick, dashed] (2, 2.7) -- (3,3);
\draw[thick, dashed] (3, 3.5) -- (3,3);
\draw[above] (3,3.5) node {$(U^n)^{+}$};
\fill (canvas cs:x=3cm,y=3.5cm) circle (1.5pt);

\fill (canvas cs:x=4.5cm,y=3.6cm) circle (1.5pt);
\draw[below] (4.5, 3.4) node {${U}^{n+1}$};
\draw[thick, dashed] (4.5, 4) -- (5.5,4.2);

\fill (canvas cs:x=4.5cm,y=4cm) circle (1.5pt);
\draw[above] (4.5, 4) node {$(U^{n+1})^{+}$};
\draw [thick, dashed] (3, 3.5) -- (4.5, 3.65);
\draw [thick, dashed] (4.5,4) -- (4.5, 3.65);

\fill (canvas cs:x=1cm,y=1.5cm) circle (1.5pt);
\draw[below] (1,1.5) node {$t_0$};
\fill (canvas cs:x=3cm,y=1.5cm) circle (1.5pt);
\draw[below] (3,1.5) node {$t_n$};
\fill (canvas cs:x=4.5cm,y=1.5cm) circle (1.5pt);
\draw[below] (4.5,1.5) node {$t_{n+1}$};

\end{tikzpicture}
\caption{An illustration of the time integration with cutoff (the dashed line).
Here, $U^n$ is the numerical solution, $({U}^{n})^{+}$ is the corrected numerical solution and
$u(t)$ is the exact solution of the IBVP starting from $u^0$ at $t=t_0$.
}
\label{f-1}
\end{figure}

We consider initial-boundary value problems of parabolic PDEs which are well posed and have
a unique, nonnegative solution. Instead of employing a maximum-principle preserving scheme,
we propose to use the cutoff method that cuts off the negative values in the computed solution
at each time step and then continues the time integration with the corrected solution;
see Fig.~\ref{f-1} for a sketch of the solution procedure.
The cutoff method shares a similar idea with many projection and correction methods.
For example, projection methods \cite{Chorin1968,Teman1968}
are an effective means of enforcing the incompressibility condition in the numerical solution of time dependent
incompressible fluid flow problems. 
Projection is also used by some researchers in the numerical solution
of Hamiltonian systems to preserve the energy; e.g., see \cite{Hai2006}.
A FE/implicit Euler discretization of a Cahn-Hilliard equation with logarithmic
nonlinearity is considered in \cite{CE92} and iterative methods to solve the resulting
nonlinear systems of equations are developed to keep the solution within its physically
relevant range.
Solution compression is used in the design of maximum-principle-satisfying high order schemes
for scalar conservation laws by Zhang and Shu \cite{ZhSh2010}.
It should be emphasized that the cutoff method provides several advantages over
many positivity-preserving or maximum-principle-preserving schemes.
Its implementation is simple and requires no significant changes in the
existing code. Moreover, the cutoff procedure does not impose any constraint on the mesh and time step
which maximum-principle-preserving schemes often impose. 
However, unlike projection methods, there do not seem to exist published theoretical or numerical
studies on the cutoff method.

The objective of this paper is to provide a theoretical and numerical study of the cutoff method
for use in the computation of nonnegative solutions of parabolic PDEs. 
We shall first prove that a broad class of finite difference methods are convergent when they are incorporated with
the cutoff method for linear parabolic PDEs. Two applications are then investigated.
The first is a linear anisotropic diffusion problems which satisfies the setting of the convergence analysis.
It is known that a conventional finite difference or finite element method with a uniform or quasi-uniform mesh
does not preserve the maximum principle (nor the solution nonnegativity) and
typically produces negative undershoot in
the computed solution. The cutoff method removes the unphysical negative values in the solution
while keeping the same convergence order of the underlying scheme.
The second application is nonlinear lubrication-type equations. It is unclear if the convergence analysis
(which can be extended to some nonlinear equations; see Remark~\ref{rem2.0})
applies to those equations. A distinct feature of lubrication-type equations is that a positive solution can
develop a finite time singularity of the form $u \to 0$ and become identically zero on some spatial interval
for a period of time.  A conventional finite difference or finite element method can
produce negative values in the computed solution during this development of singularity.
Negative solution values are not only unphysical but also cease the computation
for typical situations where the nonlinear diffusion coefficients are undefined for negative solution values.
The application of the cutoff method avoids this difficulty. Moreover, it is shown that
no regularization is needed when cutoff is used. Numerical results are
shown to be in good agreement with existing ones in the literature.

The outline of the paper is as follows. The cutoff method is described and some of its properties are given in
Section \ref{SEC:cutoff}. The convergence of a broad class of finite difference methods is proved also in
Section \ref{SEC:cutoff} when they are incorporated with the cutoff method for the computation of
nonnegative solutions of linear parabolic PDEs. Application of finite difference methods with
cutoff to a linear anisotropic diffusion problem and a lubrication-type PDE is studied in Sections~\ref{SEC:aniso-diffusion}
and \ref{SEC:lubrication}, respectively. Finally, Section~\ref{SEC:conclusion} contains conclusions and 
further comments.

\section{The cutoff method and convergence analysis}
\label{SEC:cutoff}

In this section, we first describe the cutoff method and some of its properties.
We then give a convergence analysis of the method for a class of finite difference methods
applied to linear parabolic problems. Possible generalizations are also discussed.

\subsection{The cutoff method}

Consider a continuous function $f = f(\V{x})$ defined on a domain $\Omega \subset \mathbb{R}^d$ ($d \ge 1$).
For any given cutoff parameter $\delta \in \mathbb{R}$, we define the $\delta$-cutoff function as
\beq
f_\delta^{+}(\V{x}) = 
\begin{cases}
f(\V{x}),& \text{ for } f(\V{x}) \ge \delta  \\
\delta, & \text{ for } f(\V{x}) \le \delta 
\end{cases}
\quad \forall \V{x}\,   \in \Omega .
\label{corect-2}
\eeq
Note that $f^{+}(\V{x}) \equiv f_0^{+}(\V{x})$ is the nonnegative part of $f$, i.e.,
\beq
f^{+}(\V{x}) = \frac{1}{2} (|f(\V{x})| + f(\V{x})) = 
\begin{cases}
f(\V{x}),& \text{ for } f(\V{x}) \ge 0 \\
0, & \text{ for } f(\V{x}) \le 0
\end{cases}
\quad \forall \V{x}\,   \in \Omega .
\label{corect-1}
\eeq

The following three lemmas give some properties of cutoff functions.
Although these properties are easy to prove (and almost obvious), they play a central role
in the convergence analysis for a class of finite difference methods in the next subsection.

\begin{lem}
\label{lem2.1}
For a given function $f$ and any nonnegative continuous function $u$ defined on $\Omega$, we have
\beq
|f^{+}(\V{x}) - u(\V{x})| \le |f(\V{x}) - u(\V{x})|,\quad \forall \, \V{x} \in \Omega.
\label{lem2.1-1}
\eeq
\end{lem} 

{\bf Proof.}
Obviously, (\ref{lem2.1-1}) holds when $f(\V{x}) \ge 0$. When $f(\V{x}) < 0$, we have $f^{+}(\V{x}) = 0$.
By assumption, $u(\V{x}) \ge 0$. Thus, we have
\[
| f^{+}(\V{x}) - u(\V{x}) | = | u(\V{x}) | \le  |f(\V{x})| + | u(\V{x}) | = |f(\V{x}) - u(\V{x})|.
\]
\proofend

\begin{lem}
\label{lem2.2}
For a given function $f$ and any nonnegative continuous function $u$ defined on $\Omega$, we have
\beq
|f^{+}(\V{x}) - f(\V{x})| \le |u(\V{x}) - f(\V{x}) |,\quad \forall \, \V{x} \in \Omega.
\label{lem2.2-1}
\eeq
\end{lem} 

{\bf Proof.} The proof is similar to that of Lemma~\ref{lem2.1}.
\proofend

\begin{lem}
\label{lem2.3}
Given a cutoff parameter $\delta \ge 0$, for any function $f$ and any nonnegative continuous function $u$
defined on $\Omega$, we have
\bey
&&  |f_{\delta}^{+}(\V{x}) - f^{+}(\V{x})| \le \delta, \quad \forall \, \V{x} \in \Omega
\label{lem2.3-1}
\\
&& |f_{\delta}^{+}(\V{x}) - u(\V{x})| \le |f(\V{x}) - u(\V{x})| + \delta,\quad \forall \, \V{x} \in \Omega.
\label{lem2.3-2}
\eey
\end{lem} 

{\bf Proof.} Inequality (\ref{lem2.3-1}) follows from the definitions of $f^{+}$ and $f_{\delta}^{+}$.
Inequality (\ref{lem2.3-2}) follows from the triangle inequality $|f_{\delta}^{+}(\V{x}) - u(\V{x})| 
\le |f^{+}(\V{x}) - u(\V{x})| + |f_{\delta}^{+}(\V{x}) - f^{+}(\V{x})|$ and Lemma~\ref{lem2.1}.
\proofend

\subsection{Error analysis for a linear IBVP problem}
\label{SEC:2.2}

We now use the results in the previous subsection to analyze the convergence of 
a class of finite difference schemes for a general initial-boundary value problem (IBVP) in the form
\beq
\begin{cases}
\frac{\p u}{\p t} = L(u),\quad \mbox{in} \quad \Omega\times (t_0, T] \\
u = g,\quad \mbox{on} \quad \p\Omega \\
u = u^0, \quad \mbox{on} \quad \Omega \times\{ t = t_0\}
\end{cases}
\label{ibvp-1}
\eeq
where $\Omega \subset \mathbb{R}^d$ ($d \ge 1$) is the physical domain, $L$ is a linear elliptic
differential operator and $g$ and $u^0$ are given sufficiently smooth functions. We assume that the IBVP
is well posed and admits a unique nonnegative continuous solution, $u=u(t,\V{x})$. We also assume that
$L$ and $g$ do not contain $t$ explicitly for notational simplicity.

We consider a class of finite difference schemes for (\ref{ibvp-1})  in the matrix form as
\beq
B_1 U^{n+1} = B_0 U^n + F^n,
\label{fd-1}
\eeq
where $B_0$ and $B_1$ are matrices independent of $n$ and $U^n$ is an approximation of
the exact solution at $t=t_n$, i.e., $U^n = \{ U^n_j\},\; u^n = \{u_j^n\},\; 
U^n_j  \approx u_j^n \equiv u(t_n, \V{x}_j)$ for each mesh vertex $\V{x}_j$.
We assume that (\ref{fd-1}) satisfies
\bey
&& \| B_1^{-1} \| \le K_1,
\label{fd-2}
\\
&& \| B_1^{-1}B_0 \| \le (1 + K\Delta t),
\label{fd-3}
\\
&& \frac{1}{\Delta t} \left [ B_1 u^{n+1} - B_0 u^n - F^n \right ] \longrightarrow \frac{\p u}{\p t} - L(u), \quad
	\mbox{ as} \quad \Delta t(h) \to 0 
\label{fd-4}
\eey
where $\Delta t$ and $h$ are the maximal time step size and the maximal element size, respectively,
$K$ and $K_1$ are positive constants, and $\|\cdot \|$ is a proper matrix norm.
Condition (\ref{fd-2}) requires that (\ref{fd-1}) produce a bounded solution while
(\ref{fd-3}) and (\ref{fd-4}) are the stability and consistency conditions, respectively.
The local truncation error of this scheme is defined as
\beq
{\tau}^n = B_1 {u}^{n+1} - B_0 {u}^n - F^n .
\label{lte-1}
\eeq
Assume that scheme (\ref{fd-1}) is $(p,q)$-order for some positive integers $p$ and $q$.
Then, there exists a constant $C_{lte}(u)$ (depending only on
the exact solution) such that
\beq
\| {\tau}^n \| \le C_{lte} (u) \, \Delta t \, (\Delta t^p + h^q) .
\label{lte-2}
\eeq

It is remarked that there exist many schemes satisfying assumptions (\ref{fd-2}) -- (\ref{fd-4}).
Examples include those employing central finite differences
for spatial discretization and the $\theta$-method for temporal discretization;
e.g., see Morton~and~Mayers~\cite{MM05}.

\begin{thm}
\label{thm2.1}
Assume that IBVP (\ref{ibvp-1}) is well posed and admits a unique nonnegative continuous exact solution $u=u(t,\V{x})$.
We also assume that scheme (\ref{fd-1}) satisfies (\ref{fd-2}) -- (\ref{fd-4}) and (\ref{lte-2}).
Then, the error for the cutoff solution procedure shown in Fig.~\ref{f-1} with scheme (\ref{fd-1}) is bounded by
\beq
\|(U^n)^{+} - u^n\| \le \| U^{n} - u^n \| \le \frac{K_1}{K}\, e^{K t_n} \, C_{lte}(u)\, (\Delta t^p + h^q) .
\label{thm2.1-1}
\eeq
\end{thm}

{\bf Proof.} Let $e^n = U^n - u^n$. Notice that the cutoff solution procedure shown in Fig.~\ref{f-1}
with (\ref{fd-1}) satisfies
\beq
B_1 U^{n+1} = B_0 ({U}^n)^{+} + F^n.
\label{thm2.1-2}
\eeq
Combining this with (\ref{lte-1}), we obtain
\[
B_1 e^{n+1} = B_0 (({U}^n)^{+} - u^n) - {\tau}^n .
\]
It follows that
\[
\| e^{n+1}\| \le \|B_1^{-1} B_0\| \cdot \| ({U}^n)^{+} - u^n \|  +   \| B_1^{-1}\| \cdot \| {\tau}^n\| .
\]
Combining this with (\ref{fd-2}) and (\ref{fd-3}) leads to
\beq
\| e^{n+1}\| \le (1+K\Delta t) \, \| ({U}^n)^{+} - u^n \|  +   K_1 \, \| {\tau}^n\| .
\label{thm2.1-5}
\eeq
Lemma~\ref{lem2.1} implies
\beq
 \| ({U}^n)^{+} - u^n \|  \le \| {U}^n - u^n \| .
 \label{thm2.1-7}
\eeq
Thus, we get
\beq
\| e^{n+1}\| \le (1+K\Delta t) \, \| {e}^n \|  +   K_1 \, \| {\tau}^n\| .
\label{thm2.1-6}
\eeq
Then it is standard to show that (\ref{thm2.1-1}) follows from (\ref{lte-2}) and (\ref{thm2.1-6}).
\proofend

\vspace{10pt}

\begin{rem} {\em
\label{rem2.0}
From the above proof we can see that the key condition is the convergence of
the original scheme (without cutoff). If it (without cutoff) is convergent,
using Lemma~\ref{lem2.1} we can readily show that the scheme with cutoff is also convergent.
This observation implies that the convergence analysis of the cutoff method can be extended to
more general linear or nonlinear IBVPs.
}\proofend
\end{rem}

\vspace{10pt}

When a strictly positive solution is desired, we can replace $(U^n)^{+}$ with $(U^n)^{+}_\delta$ for some positive
constant $\delta$. The following theorem can be proven
using Lemma~\ref{lem2.3} in a similar way as for Theorem~\ref{thm2.1}.

\begin{thm}
\label{thm2.2}
Assume that IBVP (\ref{ibvp-1}) is well posed and admits a unique nonnegative continuous exact solution $u=u(t,\V{x})$.
We also assume that scheme (\ref{fd-1}) satisfies (\ref{fd-2}) -- (\ref{fd-4}) and and (\ref{lte-2}).
For the cutoff solution procedure shown in Fig.~\ref{f-1} with scheme (\ref{fd-1}) and with $(U^n)^{+}$ being replaced
with $(U^n)^{+}_\delta$ for some positive $\delta$, the error is bounded by
\beq
\|(U^n)^{+}_\delta  - u^n\| \le \frac{K_1}{K}\, e^{K t_n} \, C_{lte}(u)\, (\Delta t^p + h^q)
+ \left ( \frac{(1+K \Delta t)}{K \Delta t} \, e^{K t_n} \, + 1\right ) \delta .
\label{thm2.2-1}
\eeq
\end{thm}

\begin{rem}{\em 
\label{rem2.1}
Inequality (\ref{thm2.2-1}) shows that $\delta$ should be chosen proportional to $\Delta t$
for the error bound to stay bounded as $\Delta t \to 0$. Ideally, $\delta$ should be
at the same level as the local truncation error, i.e.,
\beq
\delta = \mathcal{O}( \Delta t (\Delta t^p + h^q)) .
\label{thm2.2-2}
\eeq
This way, the terms on the right-hand side of (\ref{thm2.2-1}) have the same convergence order
as $\Delta t (h) \to 0$.
}\proofend
\end{rem}


%




\section{An anisotropic diffusion problem}
\label{SEC:aniso-diffusion}

In this section we present numerical results obtained for a 2D linear anisotropic diffusion problem
by the cutoff method described in the previous section.  The problem takes the form of 
IBVP (\ref{ibvp-1}) with
\beq
L(u) = \nabla \cdot (\mathbb{D} \, \nabla u) + f, \quad \Omega = [0,1]\times [0,1],\quad
\mathbb{D} =  \left ( \begin{array}{cc} 500.5 & 480 \\ 480 & 500.5  \end{array} \right ),
\label{ibvp-2}
\eeq
and $f$, $u^0$, and $g$ are chosen such that the exact solution of the IBVP is given by
\beq
u=\frac{1}{2} \mbox{exp}(-t) ( \tanh(-15( x-y )) +1 ).
\label{exam3.1-soln}
\eeq
The problem satisfies the maximum principle and the solution stays between 0 and 1.

It is worth pointing out that unlike lubrication-type equations which we shall consider
in the next section, for this problem the computation can continue when negative values
occur in the computed solution. From this point of view, it is unnecessary to remove negative
values at each time step. Nevertheless, this problem satisfies the setting of the convergence analysis
in the previous section and provides a good example for verifying the theory and testing
the effectiveness and accuracy of the cutoff method.

We use Cartesian grids of size $J\times J$ for the physical domain,
central finite differences for spatial discretization, and a third-order
singly diagonally implicit Runge-Kutta  (SDIRK) method \cite{Ale77,Cas79}
for temporal discretization of the underlying PDE. The discretization is standard and can be shown
to be convergent with order $(3,2)$ (third order in time and second order in space).
It is also known (and is shown below) that this standard scheme does not preserve the maximum principle
and produces spurious undershoot
and overshoot in the numerical solution. The numerical results presented below are obtained with
a fixed time step size $\Delta t = 10^{-2}$, which was found to be small enough so that the temporal
discretization error is ignorable compared to the spatial discretization error. 

\begin{figure}[thb]
\centering
\includegraphics[height=3.5in]{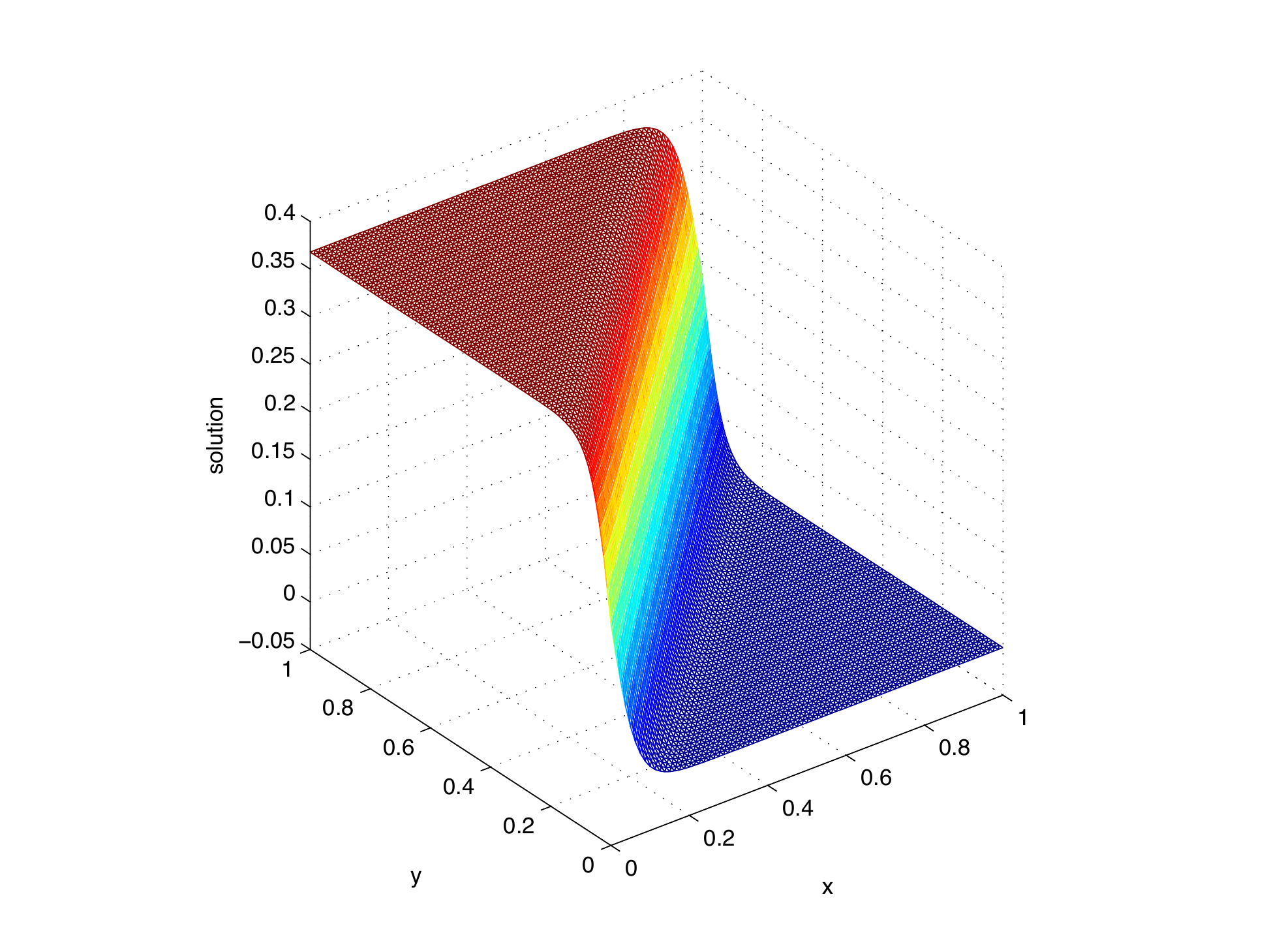}
\caption{Example~(\ref{ibvp-2}). The surface plot of a numerical solution (after cutoff) at $t=1$ obtained
with an $80 \times 80$ Cartesian grid.}
\label{Exa3.1-1}
\end{figure}

\begin{figure}[thb]
\centering
\hbox{
\begin{minipage}[b]{3in}
\centerline{(a): $U^N$}
\includegraphics[width=3.5in]{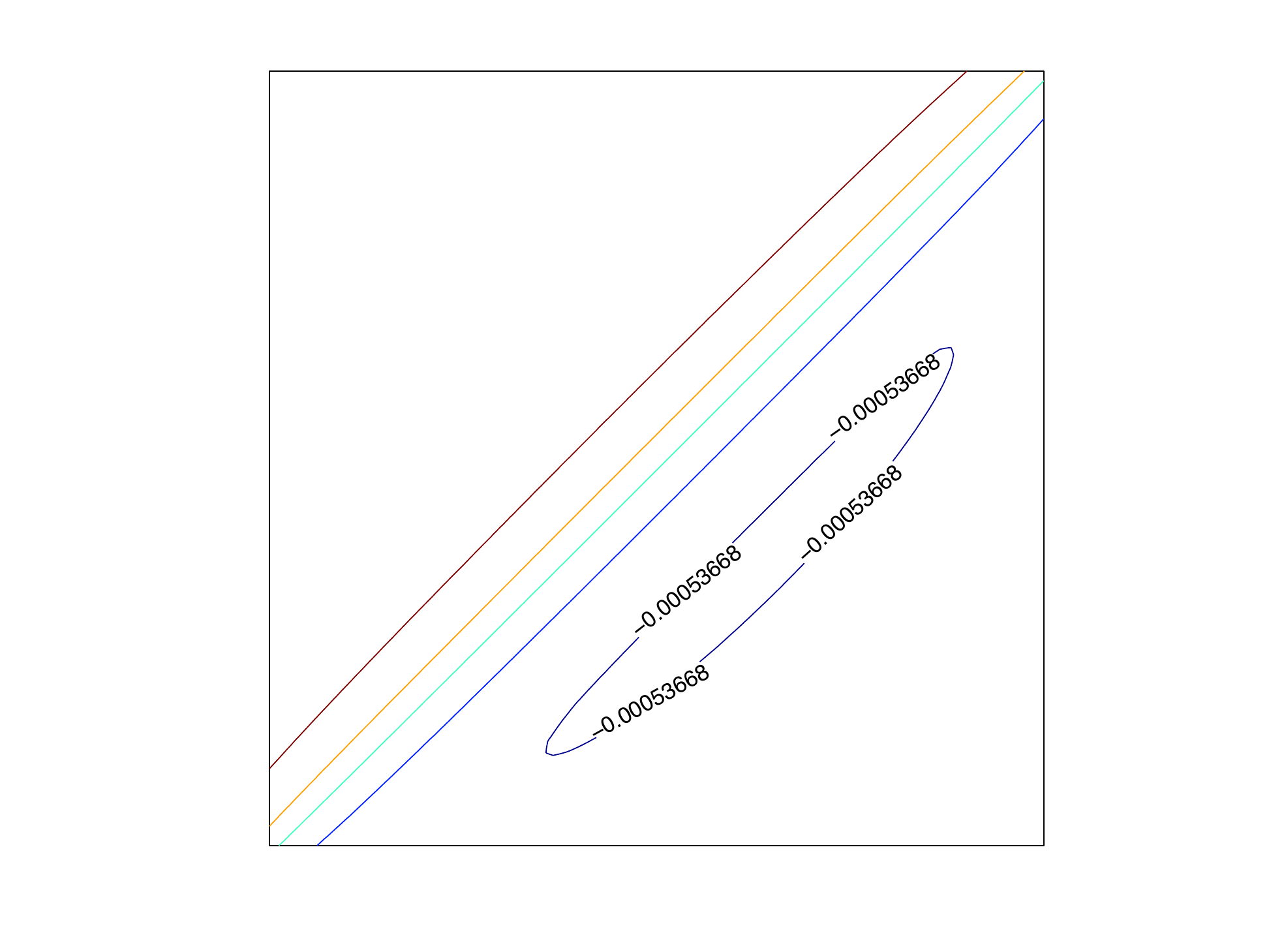}
\end{minipage}
\begin{minipage}[b]{3in}
\centerline{(b): $(U^N)^{+}$}
\includegraphics[width=3.5in]{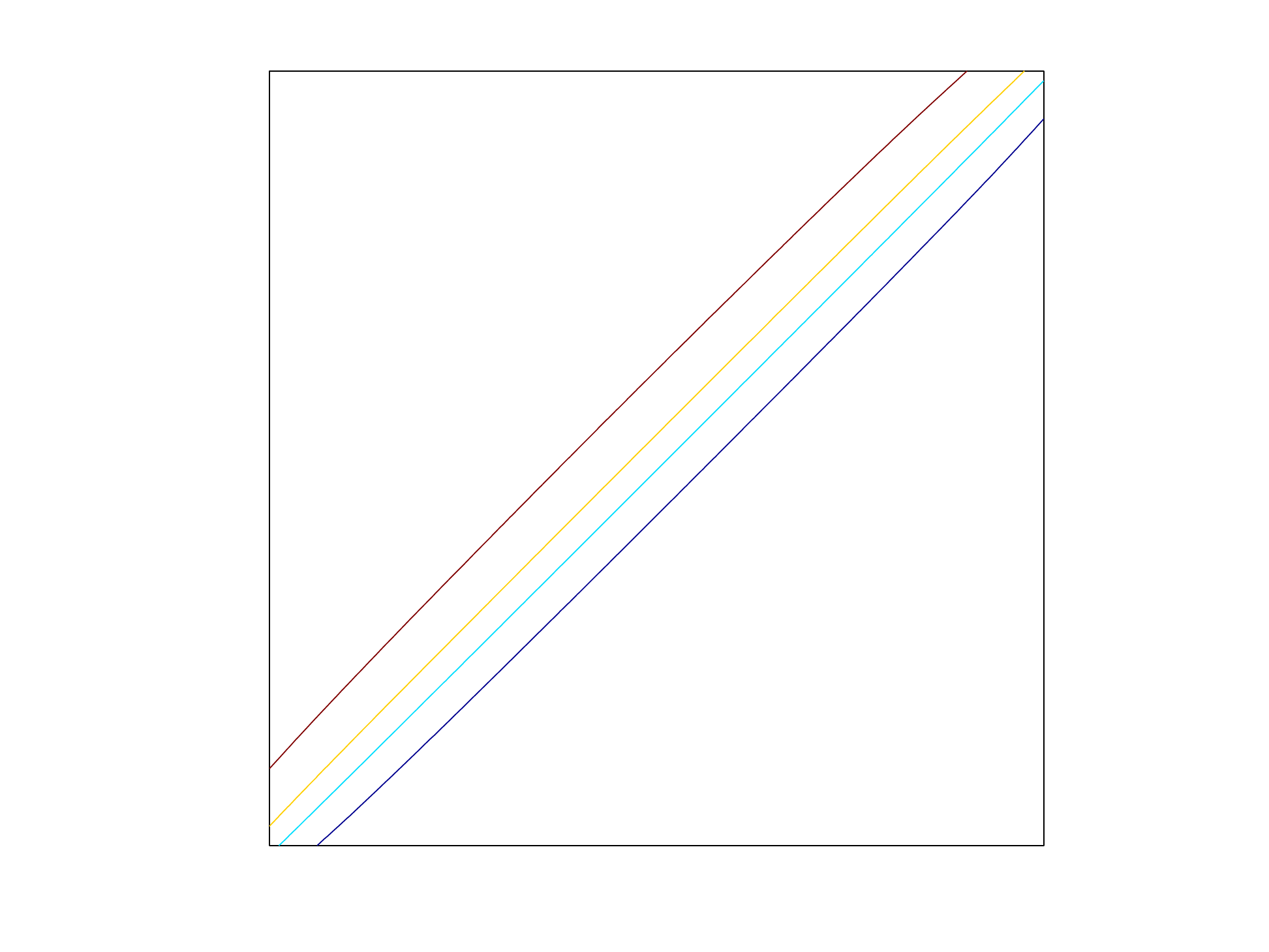}
\end{minipage}
}
\caption{Example~(\ref{ibvp-2}). Contours of the numerical solutions at $t=1$ before and after cutoff.
The solutions are obtained with an $80 \times 80$ Cartesian grid.}
\label{Exa3.1-2}
\end{figure}

Fig.~\ref{Exa3.1-1} shows the surface plot of a numerical solution (after cutoff) at $t=1$.
The contours of the numerical solutions before and after cutoff are shown in Fig.~\ref{Exa3.1-2} (a) and (b).
One may notice that the underlying finite difference scheme does not preserve the nonnegativity
and the numerical solution at $t=1$ contains negative values (with $U_{min}^N=-1.073\times 10^{-3}$
where $N$ is the last time step) before cutoff.

\begin{figure}[thb]
\centering
\includegraphics[scale=0.8]{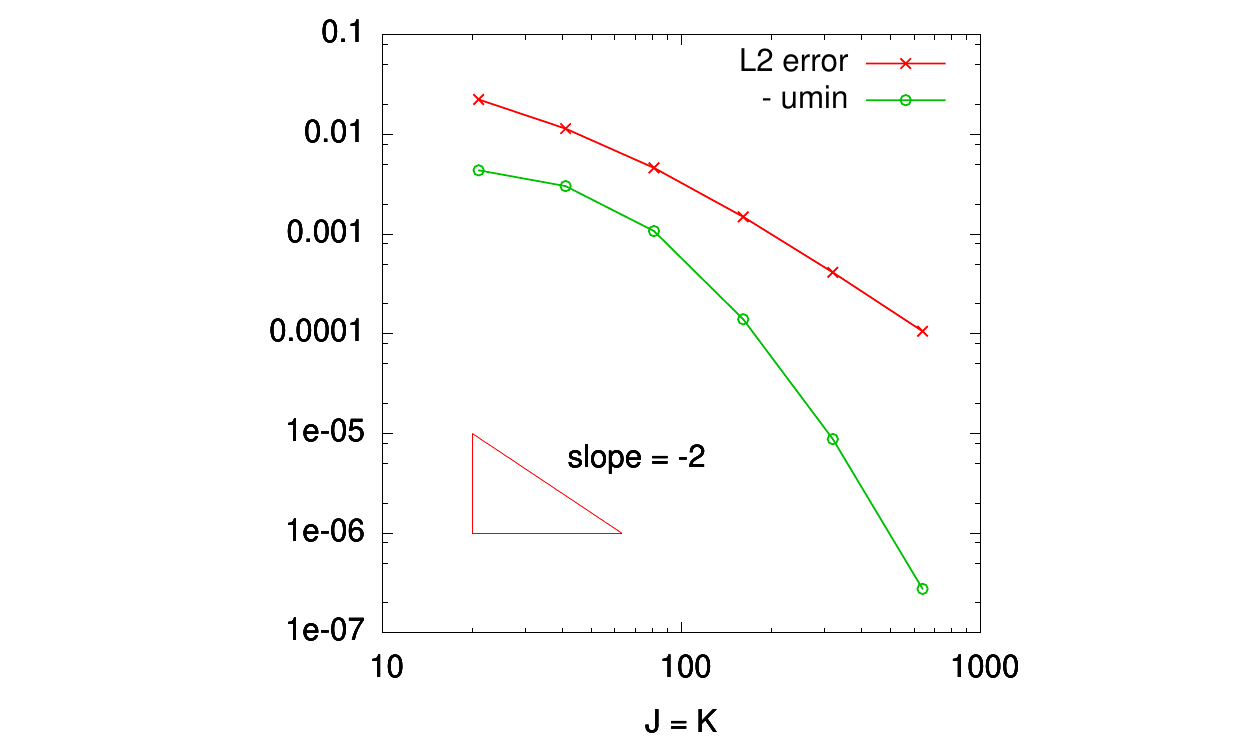}
\caption{Example~(\ref{ibvp-2}).
The $L^2$ error, $\| (U^N)^{+}-u^N \|_{L^2(\Omega)}$, and the maximal undershoot, $-u_{min}$,
are shown as functions of the number of subintervals in $x$ (or $y$) direction.
}
\label{Exa3.1-3}
\end{figure}

The $L^2$ error, $\| (U^N)^{+}-u^N \|_{L^2(\Omega)}$, and the maximal undershoot, $-u_{min}$,
are shown in Fig.~\ref{Exa3.1-3} as functions of the number of subintervals in $x$ (or $y$) direction. 
We can see that the $L^2$ error decreases at a rate of $O(J^{-2})$ (second order
in terms of element size) as $J$ increases. This is consistent with the theoretical prediction given
in Theorem~\ref{thm2.1}. On the other hand, the maximal undershoot, which is equal to the cutoff error, i.e.,
$-u_{min} = \| U^N - (U^N)^+\|_\infty$, decreases at a faster rate. This is somewhat surprising since we expect
the cutoff error to be at the level of the local truncation error, which is second order in space.

Recall that the cutoff strategy can be used to preserve the positivity of the solution. For example,
when the cutoff parameter is taken $\delta =  \Delta t ( \Delta x)^2$,
we have $(U^n_j)^{+}_{\delta} \ge \Delta t ( \Delta x)^2$. Theorem~\ref{thm2.2} implies that
the $L^2$ error, $\| (U^N)^{+}_{\delta}-u^N \|_{L^2(\Omega)}$, decreases at a rate of second order in space.
For the current example, the numerical results obtained with this cutoff parameter are indistinguishable
from those shown in Figs.~\ref{Exa3.1-2} and \ref{Exa3.1-3}. For this reason and to save space, we omit
those results here.

It is worth mentioning that computations were also performed for  
a similar problem with a convection term,
\beq
L(u) = \nabla \cdot (\mathbb{D} \, \nabla u) - \V { b } \cdot \nabla u + f, 
\label{ibvp-3}
\eeq
where $\mathbb{D}$ is given in (\ref{ibvp-2}), $\V{b} = [1000,1000]^T$, and the functions
$f$, $g$, and $u^0$ are chosen such that the exact solution is given by (\ref{exam3.1-soln}).
The obtained results (not shown) are very similar to those shown in Figs. \ref{Exa3.1-2} and \ref{Exa3.1-3}.


\section{Application to lubrication-type equations}
\label{SEC:lubrication}

In this section we consider the application of the cutoff method to a lubrication-type equation, which was first derived
from lubrication theory by Greenspan  \cite{Gre1978} for describing the movement of thin viscous films and
spreading droplets. Lubrication theory consists of a depth-averaged equation of mass conservation and a simplified
form of the Navier-Stokes equations that is appropriate for a thin layer of very viscous fluid.
We consider a lubrication-type equation in the general form
\beq
u_t + \nabla \cdot (f(u) \nabla \Delta u ) = 0,\quad f(u)\sim u^n, \;\; n \in [0, \infty )
\label{lub-1}
\eeq
where $u$ is the thickness of the viscous droplet and $n$ is a physical parameter which
has different values for different boundary conditions on the liquid solid interface.
Lubrication-type equations also appear in several other applications, including a thin neck model
in the Hele-Shaw cell \cite{CDGKS1993},  Cahn-Hilliard models  with degenerate mobility \cite{EllGar1996},
and problems in population dynamics \cite{Lew1994} and plasticity \cite{Grun1995}. 
It is noted that the solution $u$ must stay nonnegative to be physically meaningful in all of those applications.

It is known theoretically \cite{BeFr1990,BBDK1994} that in one dimension, (\ref{lub-1}) 
preserves the positivity of the solution for $n \ge 3.5$ and has a nonnegative weak solution
in a sense of distributions for $3/8 < n < 3$ (and in a weaker sense for $0 < n < 3/8$).
Such a weak solution can be obtained \cite{BeFr1990} as the limit as $\epsilon \to 0$ of the solution
of the regularized problem
\beq
u_t^\epsilon + \nabla \cdot (f^\epsilon(u^\epsilon) \nabla \Delta u^\epsilon ) = 0,\quad
f^\epsilon(u^\epsilon) = \frac{(u^\epsilon)^4 f(u^\epsilon)}{\epsilon f(u^\epsilon) + (u^\epsilon)^4} ,
\label{lub-2}
\eeq
where $\epsilon > 0$ is the regularization parameter. Moreover, numerical computation (such as see
\cite{Ber1995, Ber1996, BBDK1994}) shows that for small values of $n > 0$, a positive solution
can first develop a finite time singularity of the form $u \to 0$ at some point (which physically describes
the rupture of the liquid film), then becomes identically zero on an interval of time dependent length for a period of time,
becomes positive again at a later time, and eventually decays to the mean of the initial solution.
It is challenging to simulate this singularity development
since conventional numerical methods do not preserve the nonnegativity of the solution and requires
a huge number of mesh points to provide the necessary resolution to avoid spurious, negative undershoot.
Moreover, when $n$ is an even denominator rational number,  $f(u)$ is undefined for negative $u$.
In this situation, the computation typically stops around the initial formation of the singularity when negative values
start to appear in the numerical solution. For this reason, the simulation of the singularity development in (\ref{lub-1})
is a good test for the cutoff method although it is unclear if the convergence analysis described
in \S~\ref{SEC:2.2} applies to (\ref{lub-1}).

Great effort has been made in the numerical solution of (\ref{lub-1}).
In addition to the above mentioned references \cite{Ber1995, Ber1996, BBDK1994}, 
Zhornitskaya and Bertozzi \cite{ZhoBer2000} propose a positivity preserving finite difference scheme
for the regularized equation (\ref{lub-2}). Russell et al. \cite{RuWiXu2007} solve the same regularized
equation using a moving collocation method.
Sun et al. \cite{SunRusXu2007} solve (\ref{lub-1}) in two dimensions using an adaptive finite element method
and show that proper mesh adaptation can provide accurate resolution and there is no need 
to use regularization in the numerical simulation of the singularity development in (\ref{lub-1}).

We first consider an IBVP of the one-dimensional lubrication-type equation as
\beq
\begin{cases}
& u_t +  (f(u) u_{xxx})_x = 0, \;\; x \in \Omega \equiv (-1,1),\;\; f(u) = u^{\frac{1}{2}},
\\
& u(t,\pm 1) = 0, \; u_{xxx} (t, \pm 1) = 0,
\\
& u(0,x) = 0.8 -\cos (\pi x) + 0.25 \cos (2\pi x). 
\end{cases}
\label{ibvp-4}
\eeq
Notice that the initial condition is strictly positive (with minimum value 0.05 and mean value 0.8)
and $f(u)$ is undefined for negative $u$.
We use a uniform mesh of $J$ cells for $\Omega$, central finite
differences for spatial discretization, a third-order SDIRK method for the temporal discretization
of the underlying PDE. (A small fixed time step $\Delta t = 10^{-6}$ is used in our computation.)
This discretization is basically the same as that for the example in the previous section except
that a special treatment is needed for the diffusion coefficient for the current problem
since it is nonlinear. Recall that the scheme does not preserve the nonnegativity of the solution.
Thus, iterates can have negative values at some point during the Newton iterative solution
of the nonlinear algebraic equations resulting from the SDIRK/finite difference discretization 
of the underlying PDE.
Once this occurs, the computation stops because $f(u)$ is undefined for those values.
One way to avoid this difficulty is to use $f((U^{n+1})^{+})$ instead of $f(U^{n+1})$ in the scheme.
However, the nonsmoothness of $(U^{n+1})^{+}$ can cause difficulty in the 
convergence of the Newton iteration.
We use here the ``lagged diffusivity" method, i.e., the diffusion coefficient is computed using
the (corrected) numerical solution at the previous time step. An iterative method based on
the lagged diffusivity is used and proved to be convergent by Dobson and Vogel \cite{DoVo1997}
for a total variation denoising problem which also has a nonlinear diffusion coefficient.



\begin{figure}[thb]
\centering
\includegraphics[width=4in]{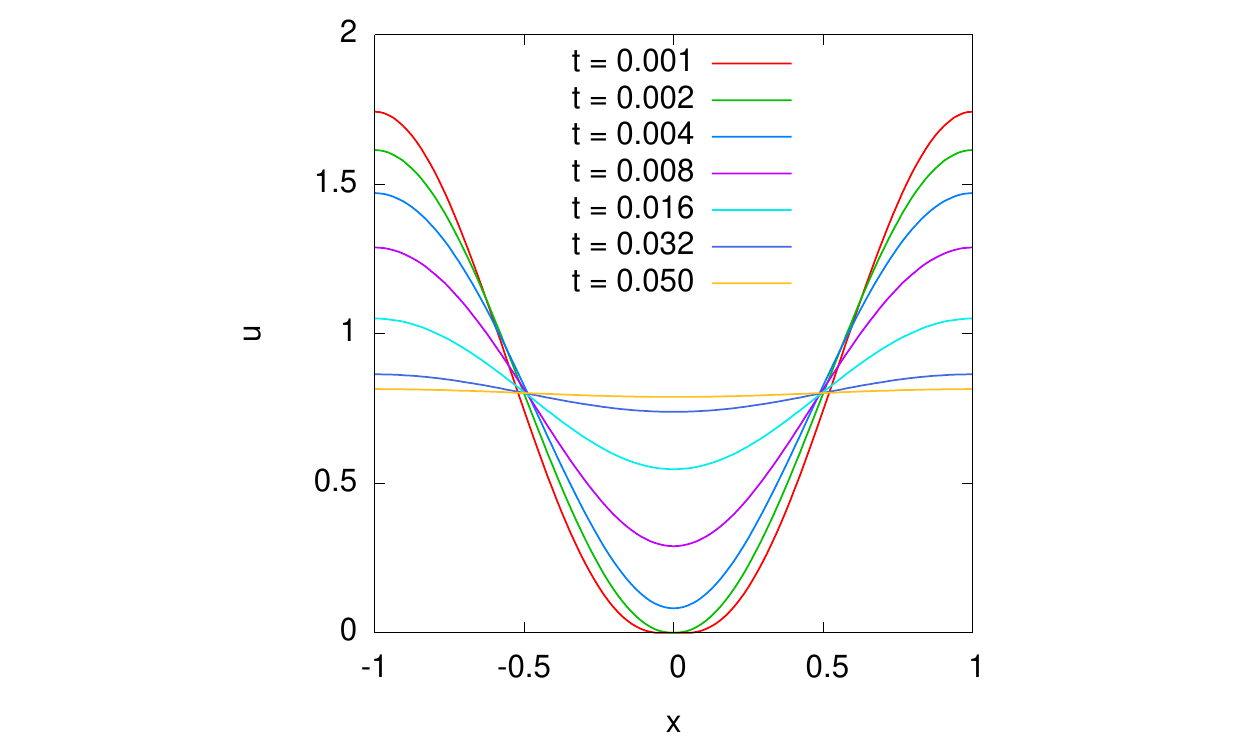}
\caption{The one-dimensional lubrication-type problem. The solution obtained with 129 uniform grid points is
shown at various time instants.}
\label{Exa4.1-1}
\end{figure}

Fig.~\ref{Exa4.1-1} shows a numerical solution at various time instants obtained with the cutoff method
(without using regularization for $f(u)$) on a relatively coarse uniform mesh of 129 points. The result is almost identical
(in the eyeball norm) to those in \cite{Ber1995,SunRusXu2007}. To further verify the result, the computation
is done with a uniform mesh of 1001 points and the obtained solution is indistinguishable from the one shown
in Fig.~\ref{Exa4.1-1} in the plotting precision.

\begin{figure}[thb]
\centering
\hbox{
\begin{minipage}[b]{2.0in}
\centerline{(a): t=7.29e-4}
\includegraphics[width=2in]{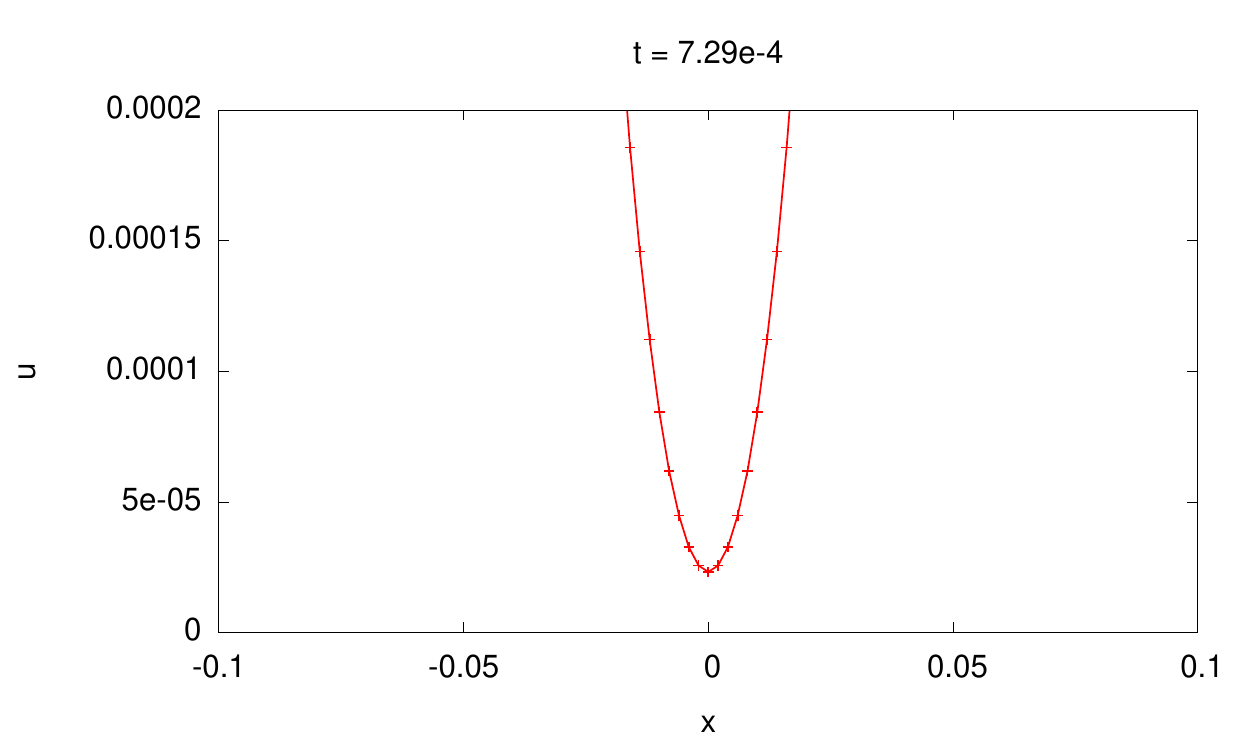}
\end{minipage}
\begin{minipage}[b]{2.0in}
\centerline{(b): t=7.5e-4}
\includegraphics[width=2in]{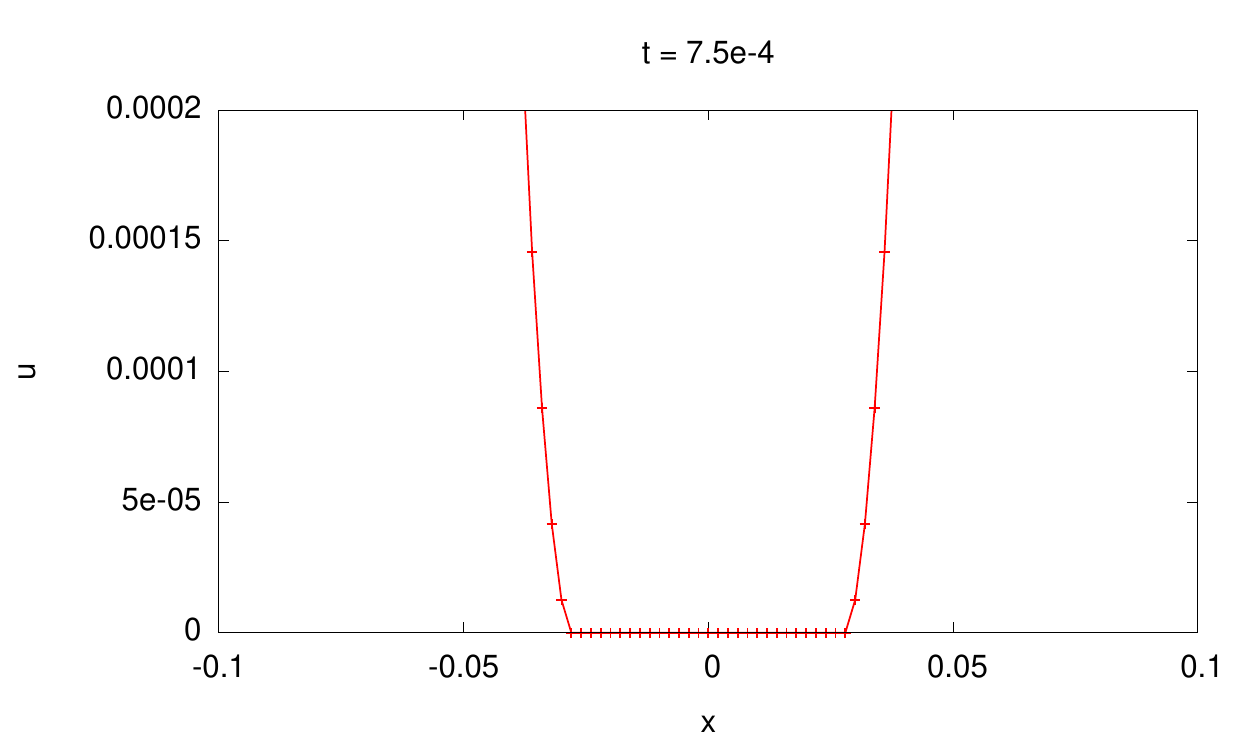}
\end{minipage}
\begin{minipage}[b]{2.0in}
\centerline{(c): t=8.0e-4}
\includegraphics[width=2in]{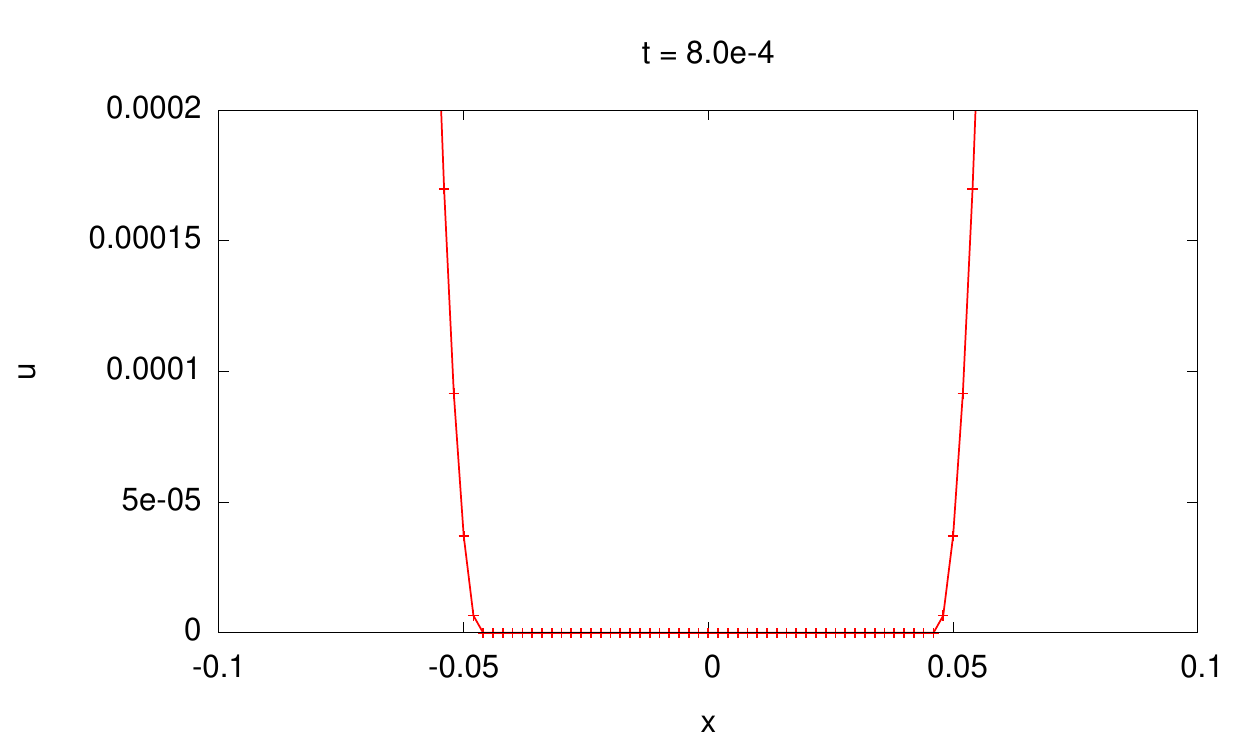}
\end{minipage}
}
\hbox{
\begin{minipage}[b]{2.0in}
\centerline{(d): t=8.5e-4}
\includegraphics[width=2in]{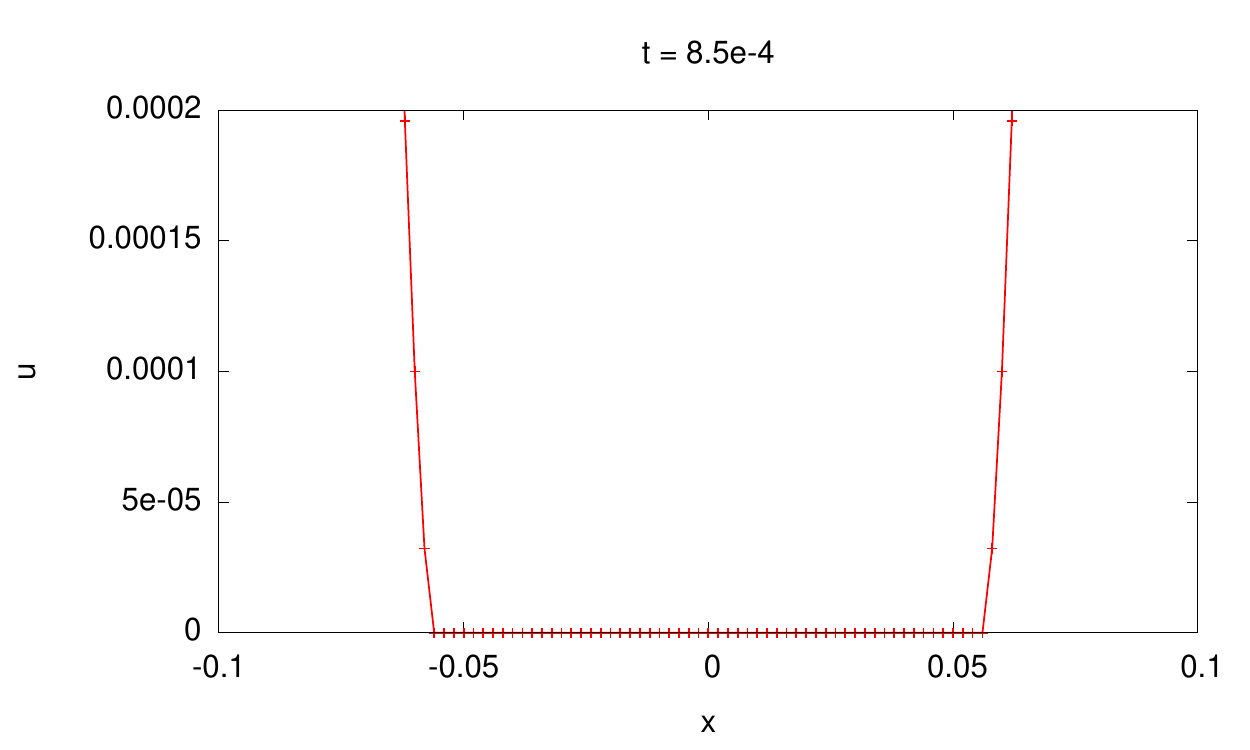}
\end{minipage}
\begin{minipage}[b]{2.0in}
\centerline{(e): t=1.0e-3}
\includegraphics[width=2in]{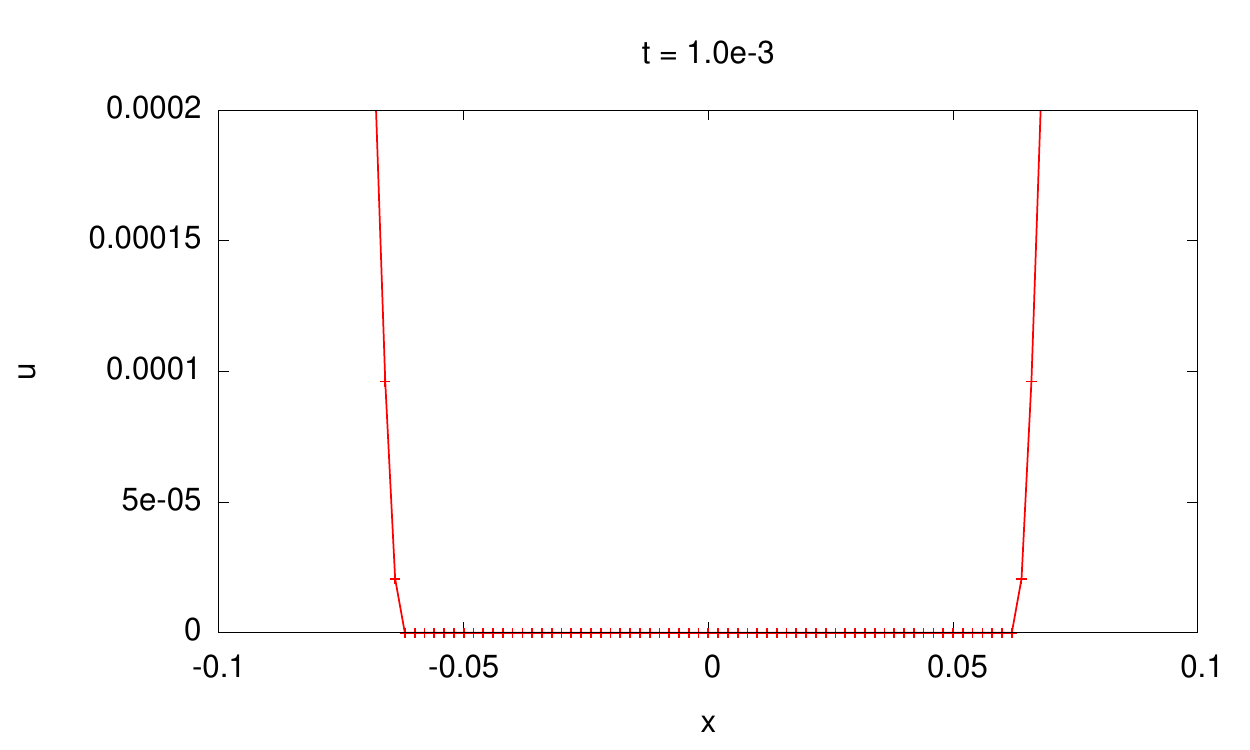}
\end{minipage}
\begin{minipage}[b]{2.0in}
\centerline{(f): t=1.5e-3}
\includegraphics[width=2in]{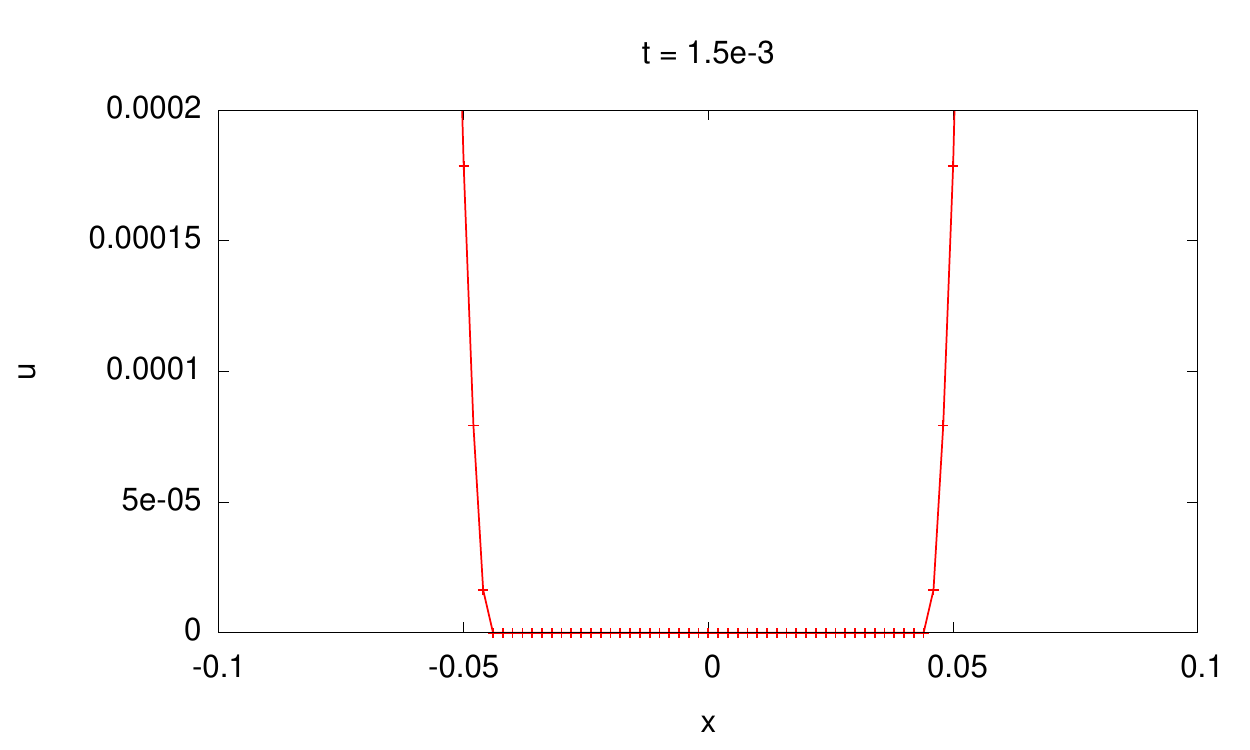}
\end{minipage}
}
\hbox{
\begin{minipage}[b]{2.0in}
\centerline{(g): t=2.0e-3}
\includegraphics[width=2in]{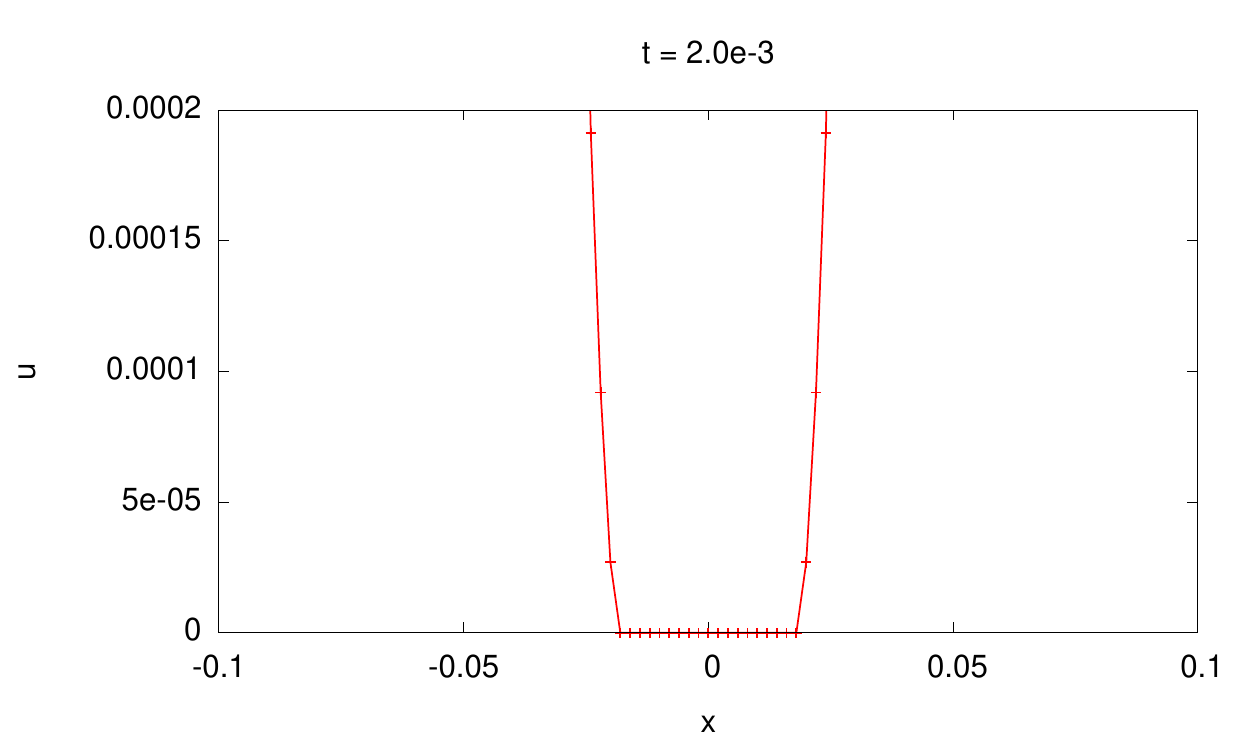}
\end{minipage}
\begin{minipage}[b]{2.0in}
\centerline{(h): t=2.3e-3}
\includegraphics[width=2in]{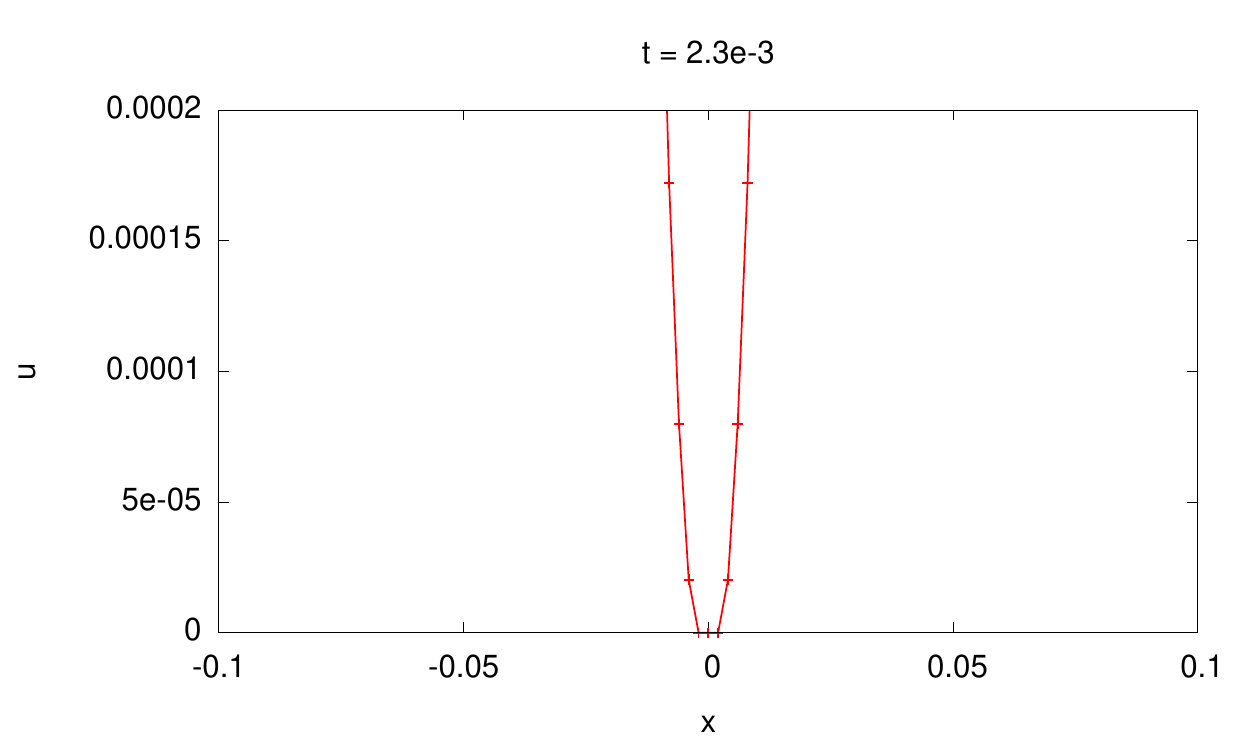}
\end{minipage}
\begin{minipage}[b]{2.0in}
\centerline{(i): t=2.34e-3}
\includegraphics[width=2in]{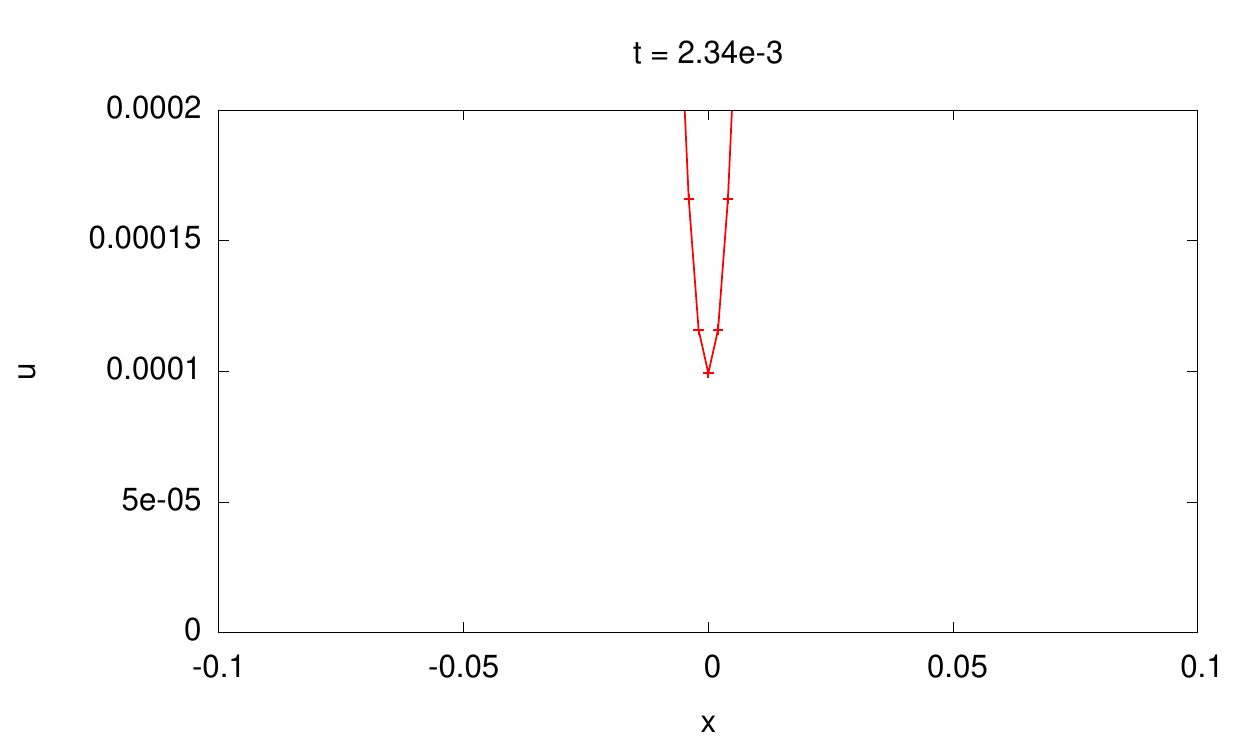}
\end{minipage}
}
\caption{The one-dimensional lubrication-type problem. The close views of the numerical solution
are shown at various time instants during the singularity development. The numerical solution is obtained
by the cutoff method (without using regularization for $f(u)$) on a uniform mesh of 1001 points
with $\Delta t = 10^{-6}$.}
\label{Exa4.1-2}
\end{figure}

\begin{figure}[thb]
\centering
\hbox{
\begin{minipage}[b]{2.0in}
\centerline{(a): t=7.29e-4}
\includegraphics[width=2in]{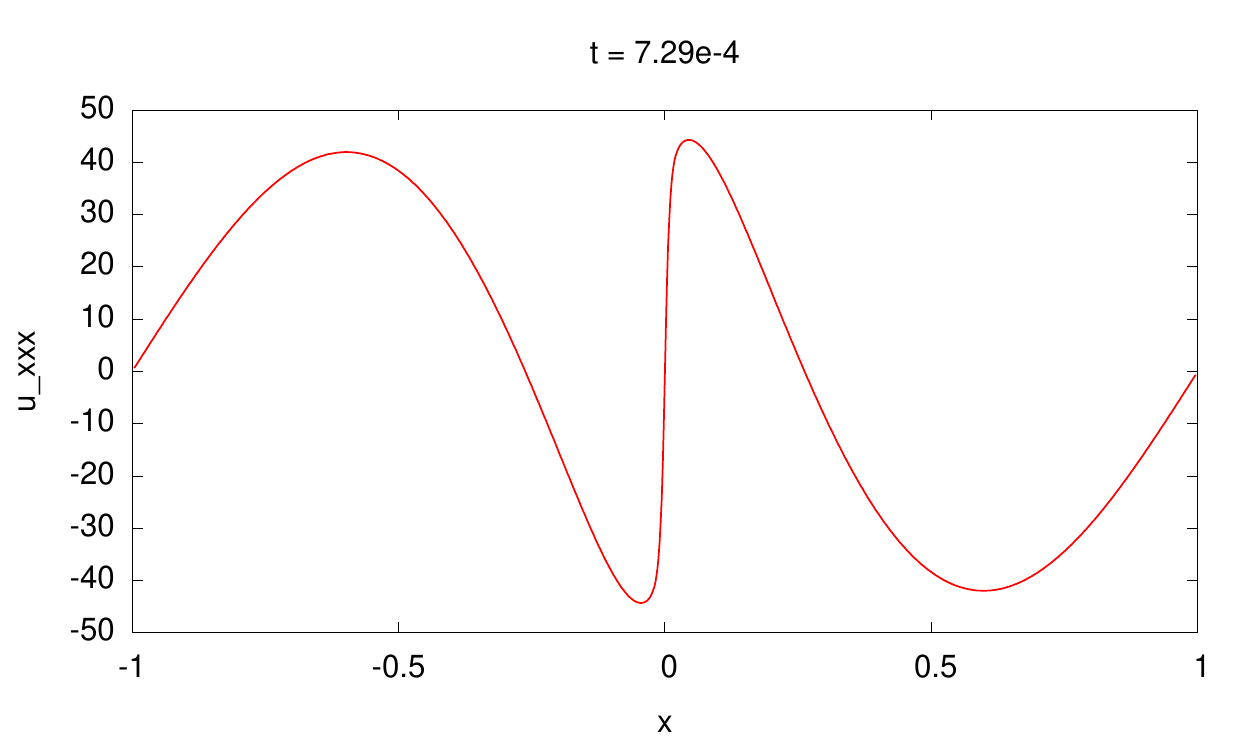}
\end{minipage}
\begin{minipage}[b]{2.0in}
\centerline{(b): t=7.5e-4}
\includegraphics[width=2in]{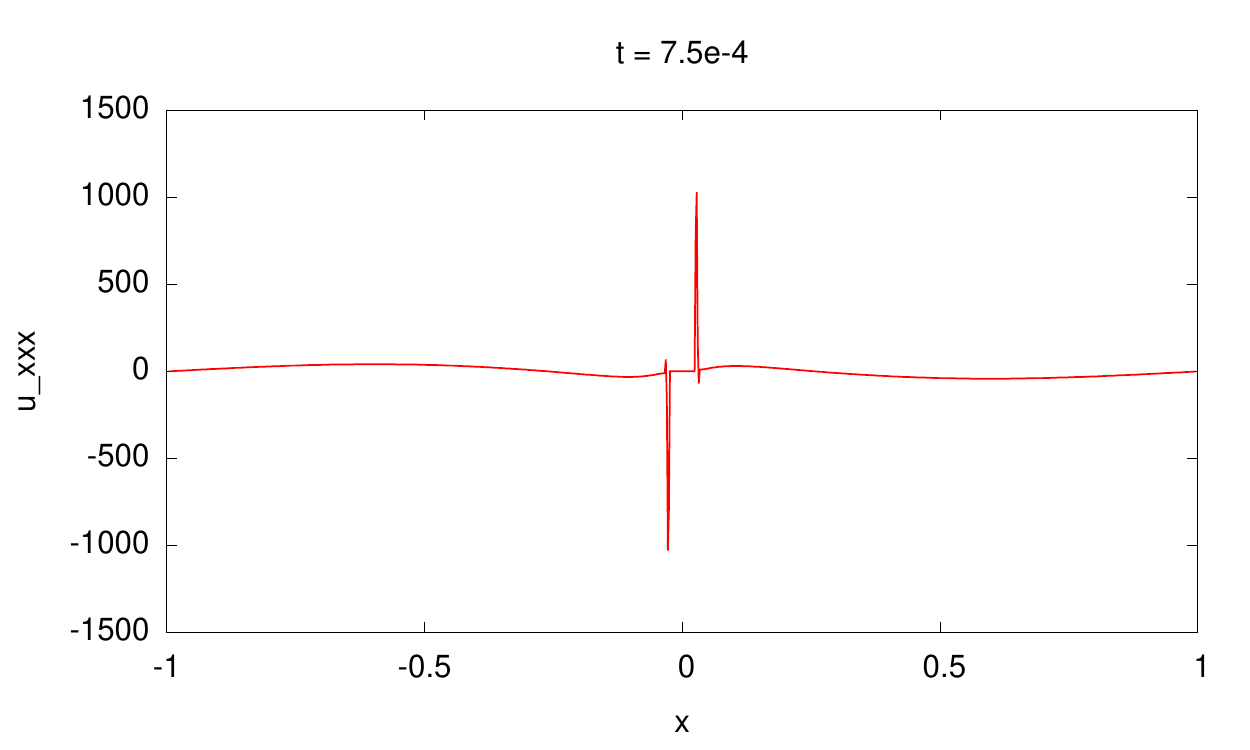}
\end{minipage}
\begin{minipage}[b]{2.0in}
\centerline{(c): t=8.0e-4}
\includegraphics[width=2in]{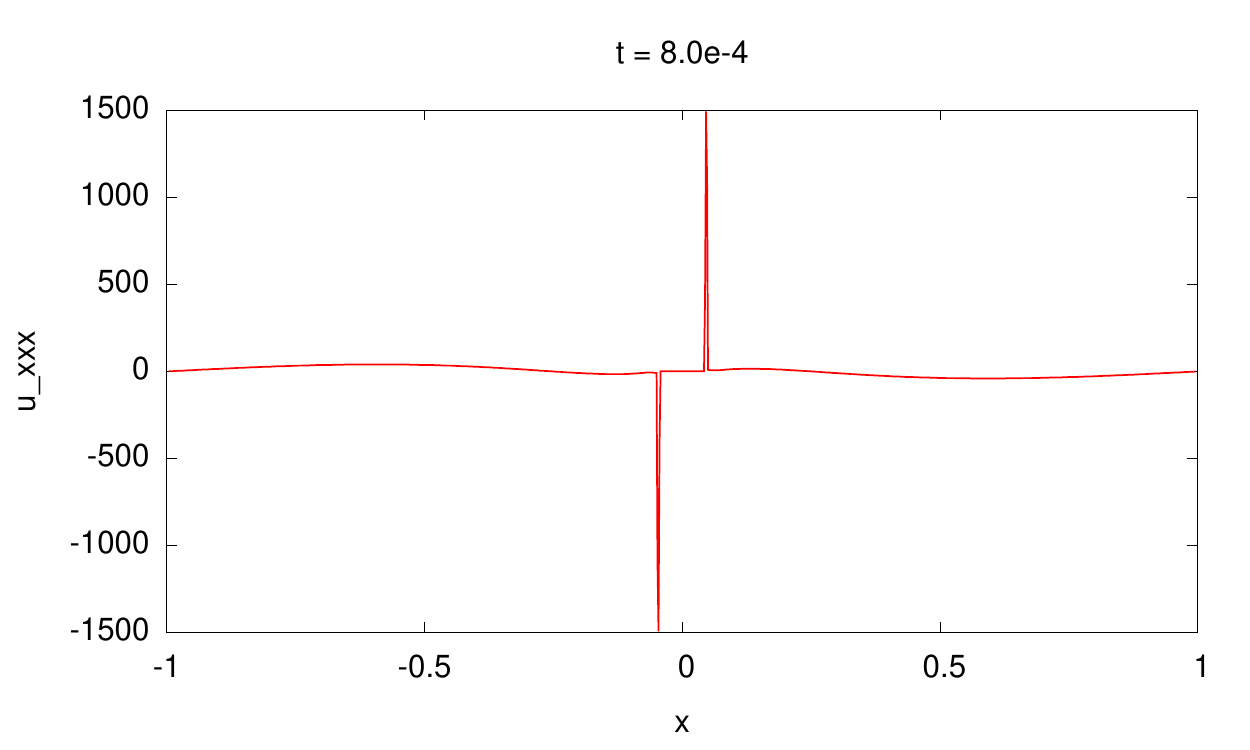}
\end{minipage}
}
\hbox{
\begin{minipage}[b]{2.0in}
\centerline{(d): t=8.5e-4}
\includegraphics[width=2in]{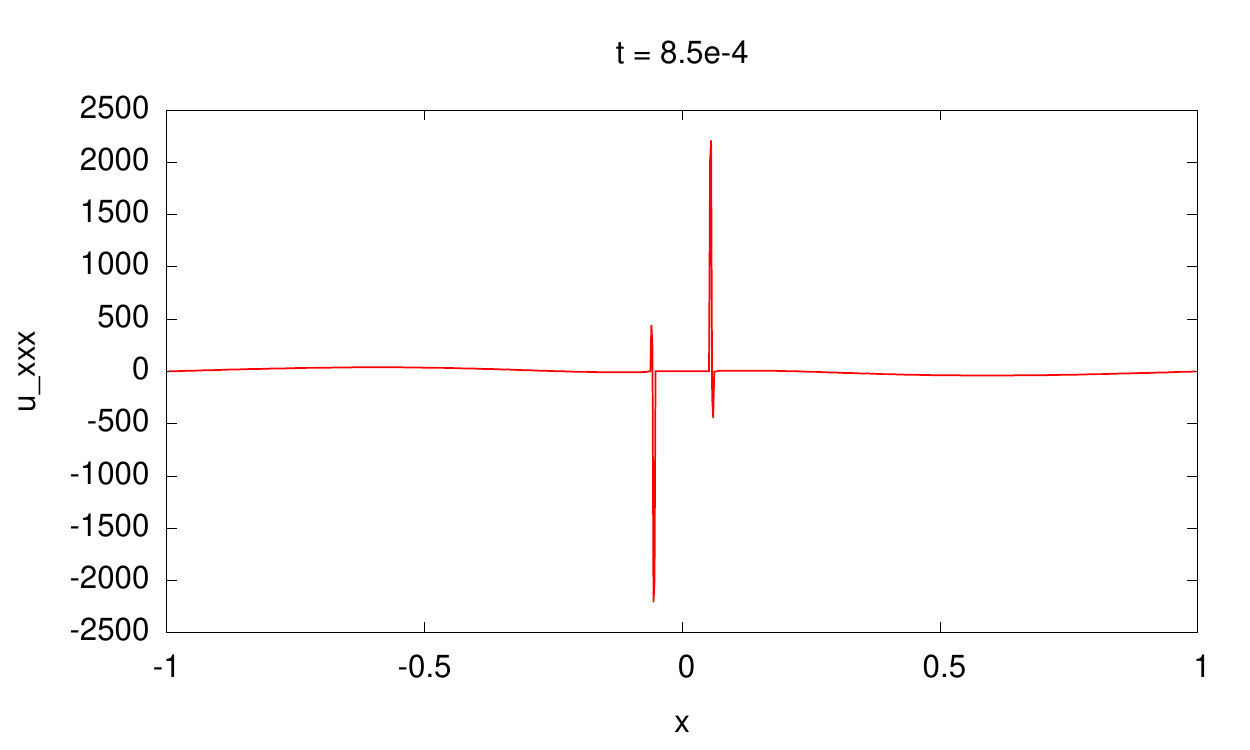}
\end{minipage}
\begin{minipage}[b]{2.0in}
\centerline{(e): t=1.0e-3}
\includegraphics[width=2in]{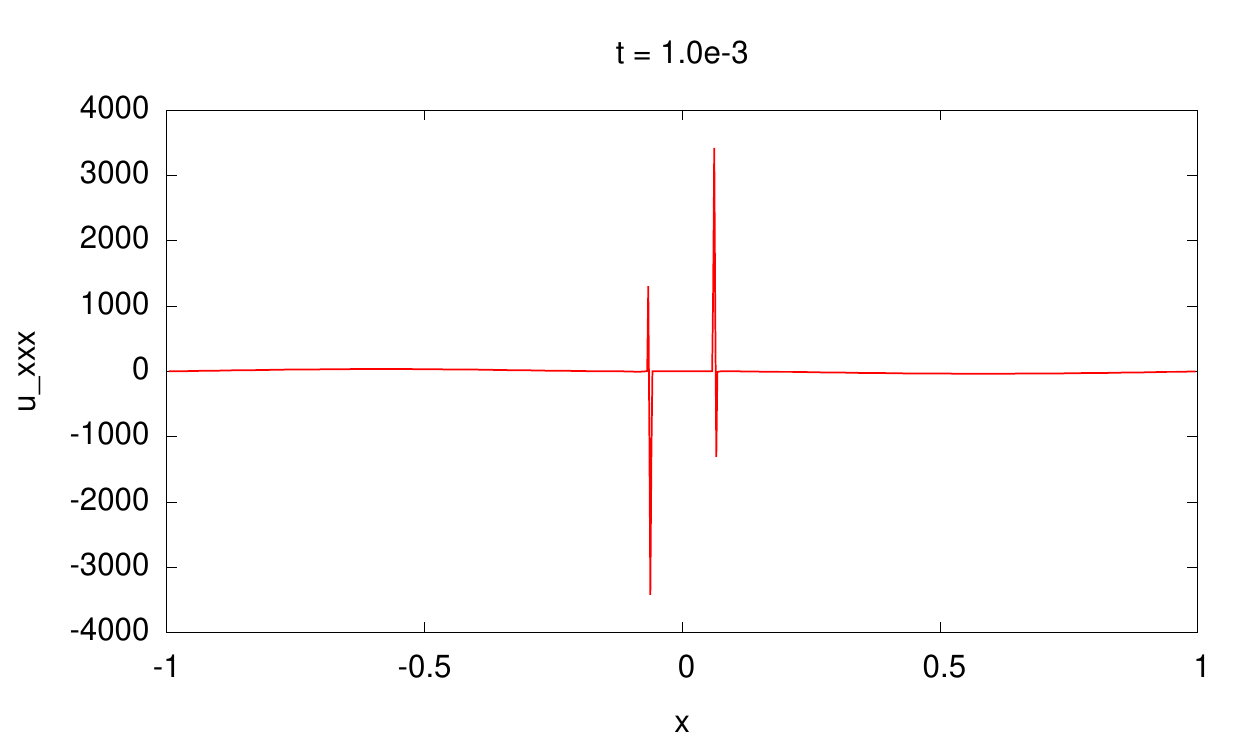}
\end{minipage}
\begin{minipage}[b]{2.0in}
\centerline{(f): t=1.5e-3}
\includegraphics[width=2in]{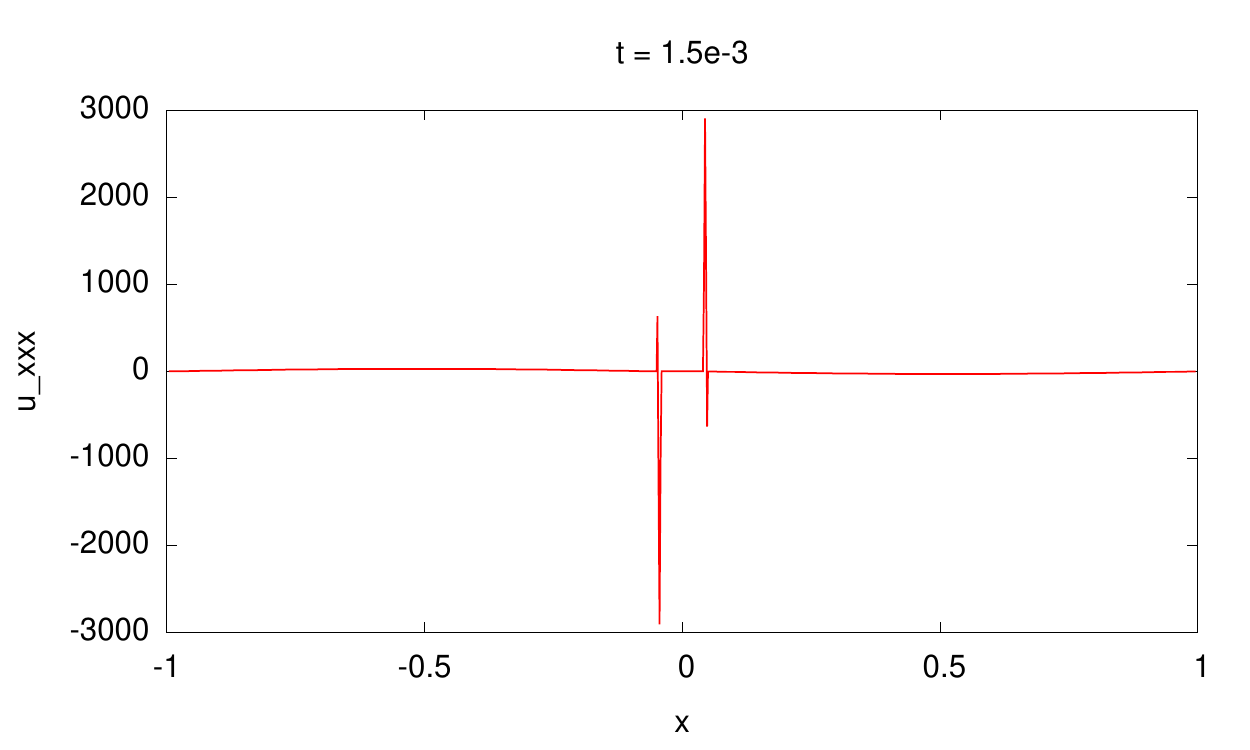}
\end{minipage}
}
\hbox{
\begin{minipage}[b]{2.0in}
\centerline{(g): t=2.0e-3}
\includegraphics[width=2in]{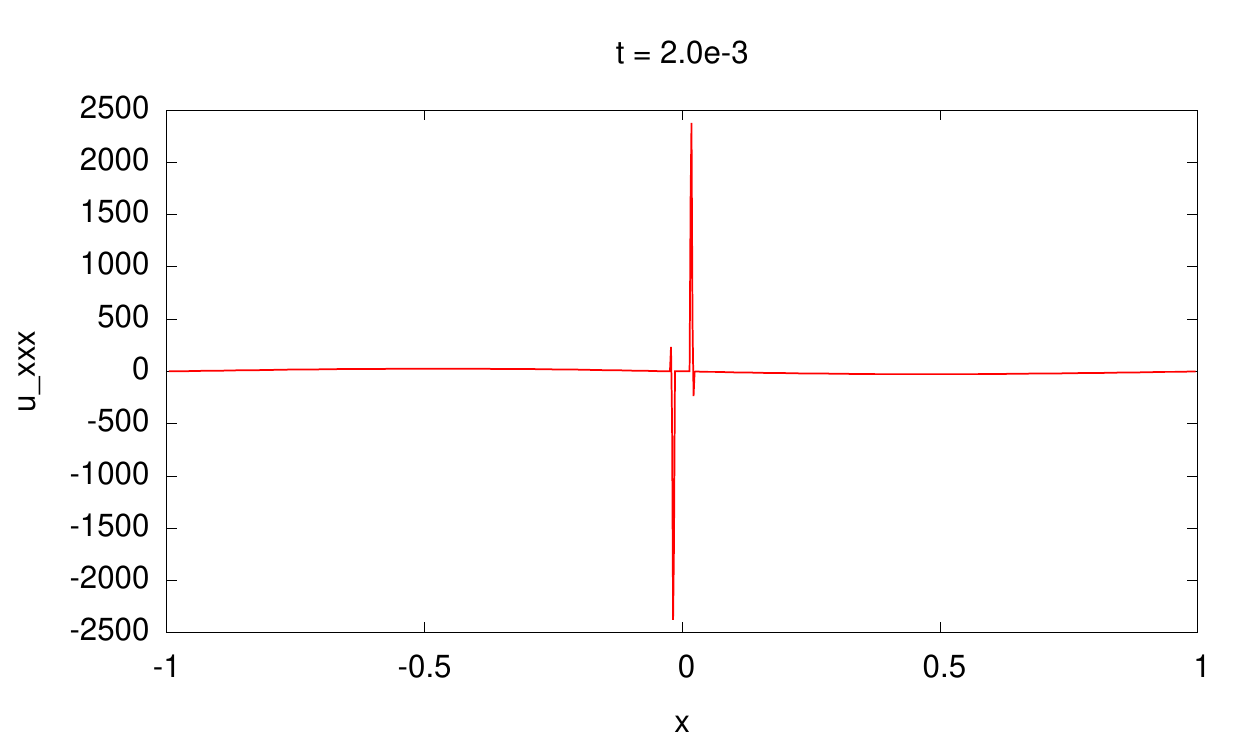}
\end{minipage}
\begin{minipage}[b]{2.0in}
\centerline{(h): t=2.3e-3}
\includegraphics[width=2in]{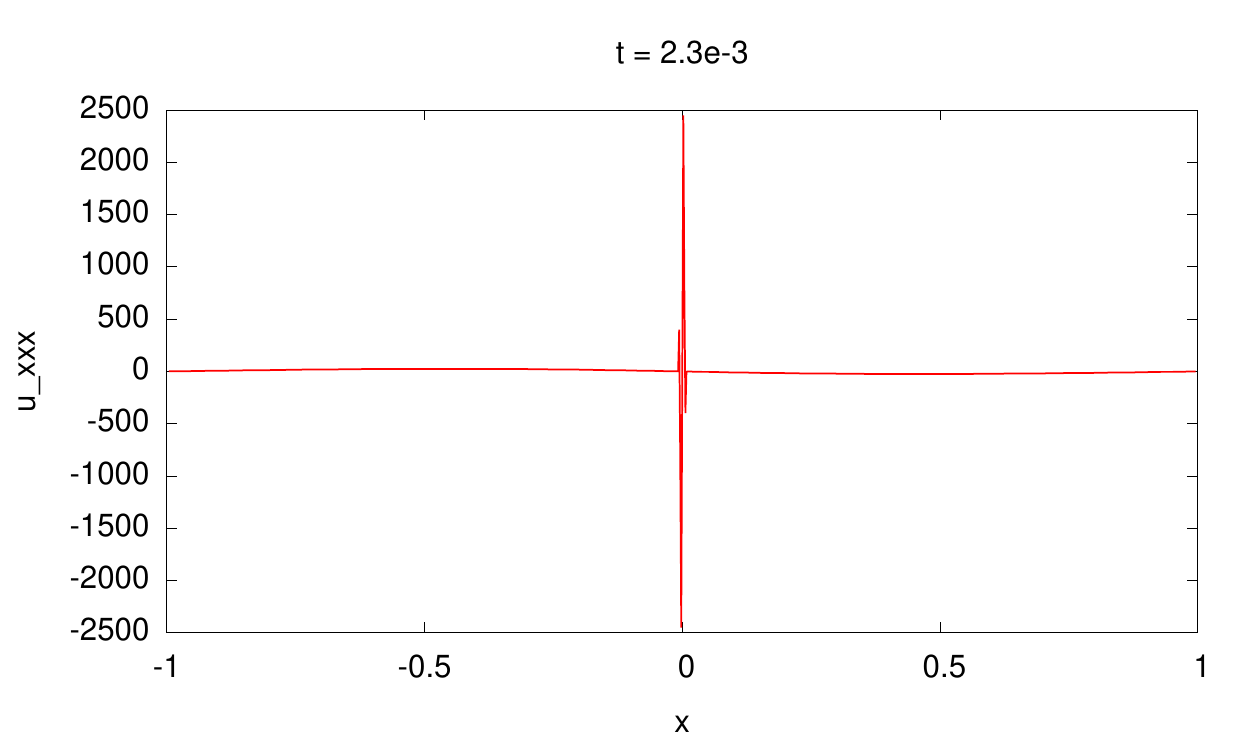}
\end{minipage}
\begin{minipage}[b]{2.0in}
\centerline{(i): t=2.34e-3}
\includegraphics[width=2in]{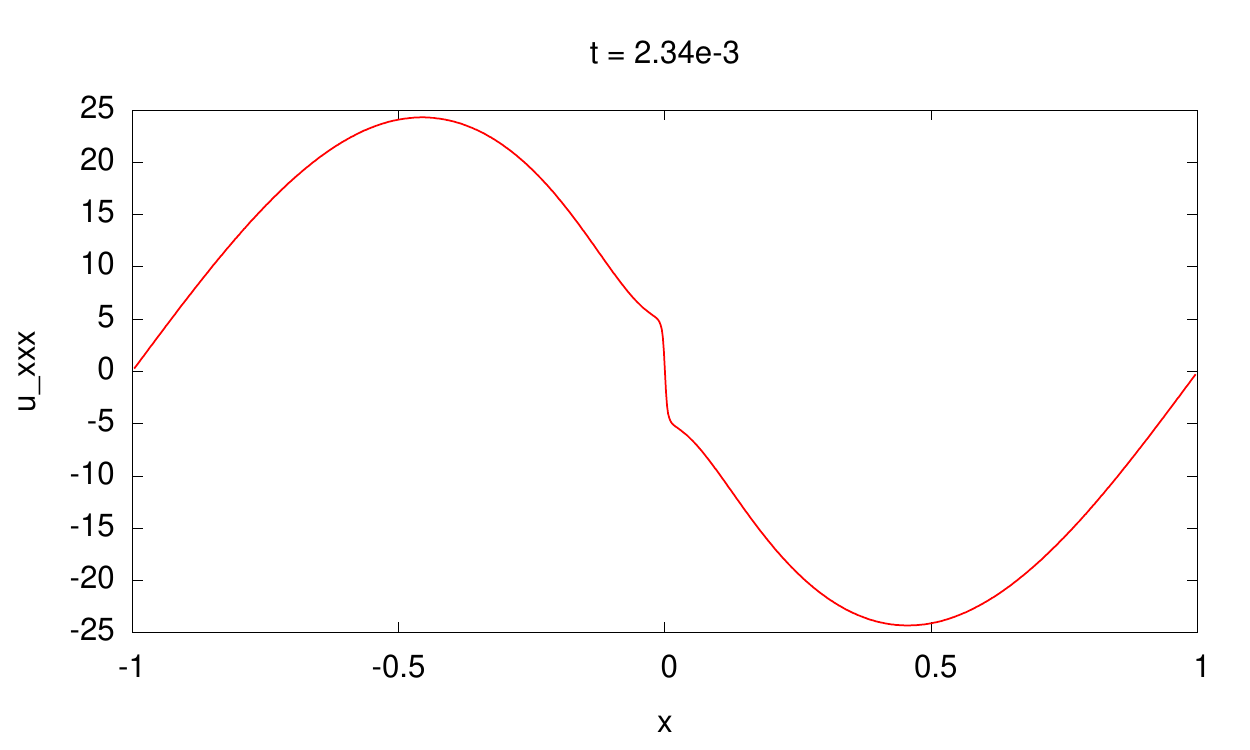}
\end{minipage}
}
\caption{The one-dimensional lubrication-type problem. The close views of the third order derivative
of the numerical solution are shown at various time instants during the singularity development.
The numerical solution is obtained by the cutoff method (without using regularization for $f(u)$)
on a uniform mesh of 1001 points with $\Delta t = 10^{-6}$.}
\label{Exa4.1-3}
\end{figure}

A uniform mesh of 1001 points is used for the numerical study of the singularity development of the solution.
Fig.~\ref{Exa4.1-2} shows the close views of the solution during the development. It can be seen that
the singularity is developed at around $t= 7.30\times 10^{-4}$, which is consistent with existing
numerical simulations,
including those in \cite{ZhoBer2000} with a positivity-preserving scheme.
Then, the solution becomes identically zero on an interval of time dependent length for a period
of time. The length of the interval increases from zero, attains its maximum (about $0.12$), and decreases to
zero. Afterward, the solution becomes positive again, and eventually decays to the mean of the initial solution.
It is known (e.g., see \cite{Ber1995}) that the 
development of the singularity is characterized by 
the third order derivative of the solution becoming infinite at certain points. 
This can be seen in Fig.~\ref{Exa4.1-3} where the third order derivative of the numerical solution is shown
at various time instants. In particular, the third order derivative becomes discontinuous when 
the singularity
occurs, then has jump discontinuities, and then becomes continuous again when the solution is positive.

\begin{figure}[thb]
\centering
\hbox{
\begin{minipage}[b]{2.0in}
\centerline{(a): t=7.3e-4}
\includegraphics[width=2in]{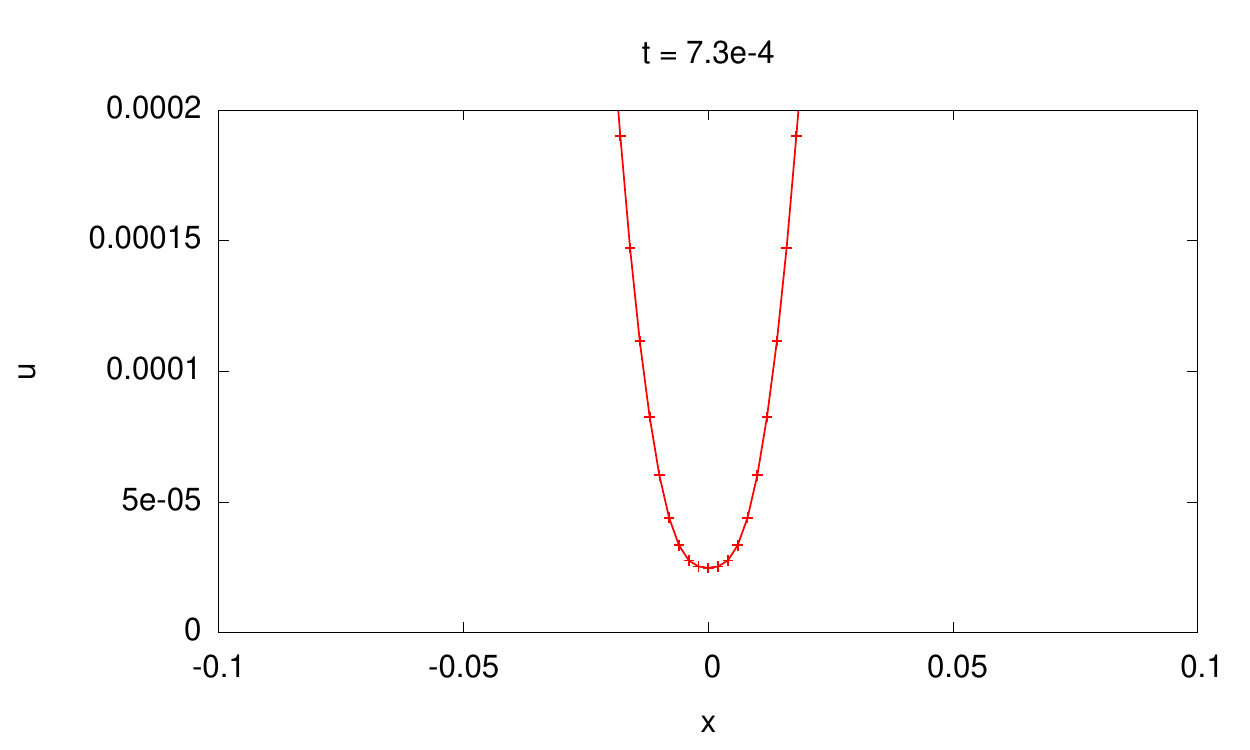}
\end{minipage}
\begin{minipage}[b]{2.0in}
\centerline{(b): t=7.5e-4}
\includegraphics[width=2in]{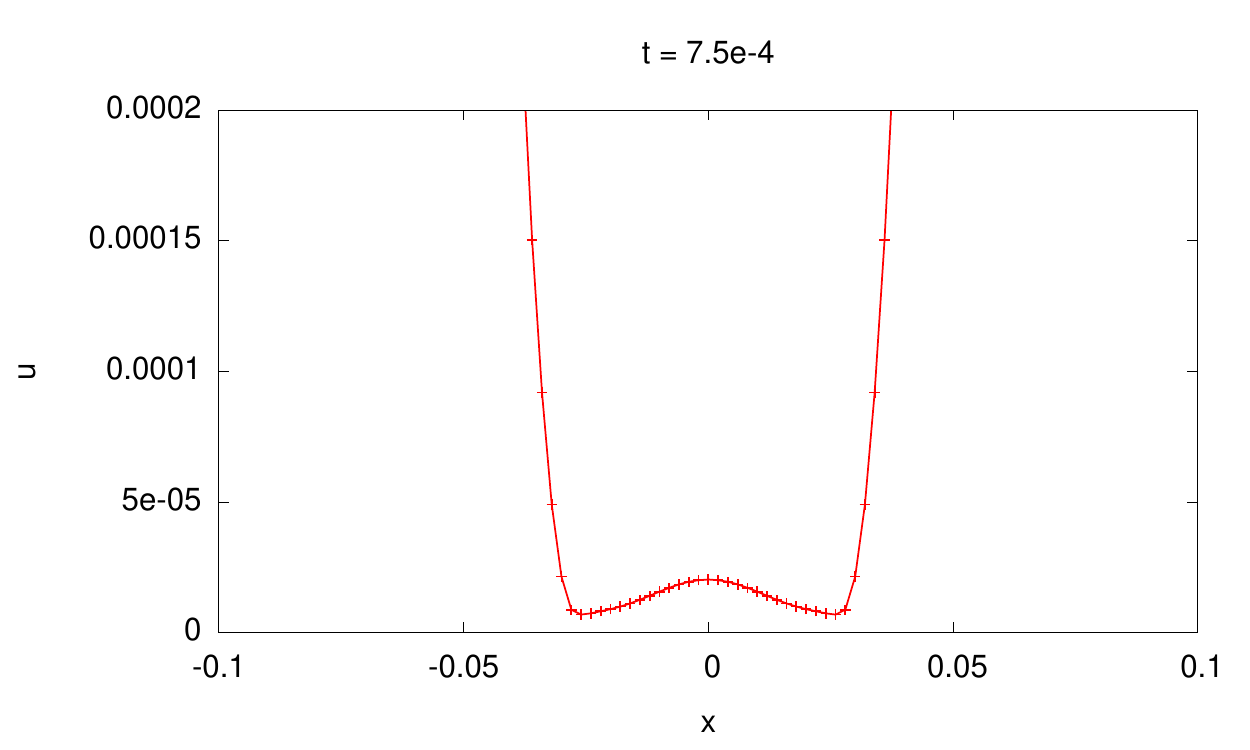}
\end{minipage}
\begin{minipage}[b]{2.0in}
\centerline{(c): t=8.0e-4}
\includegraphics[width=2in]{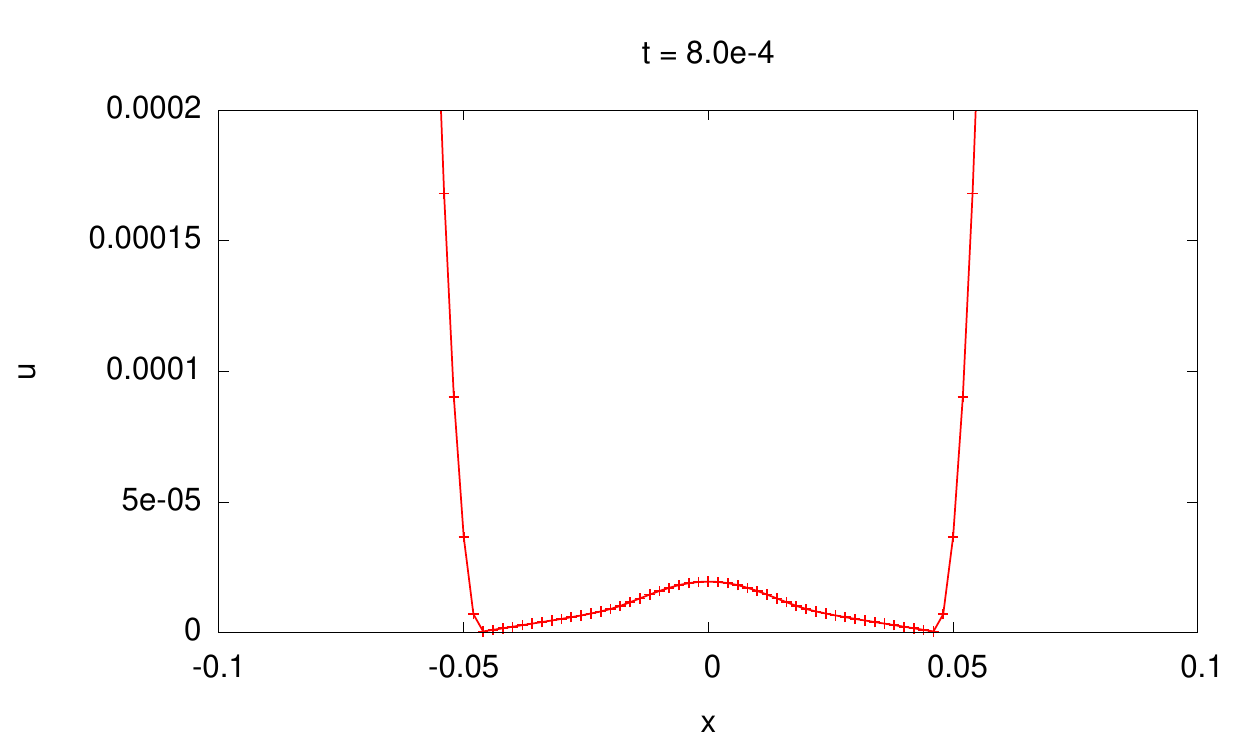}
\end{minipage}
}
\hbox{
\begin{minipage}[b]{2.0in}
\centerline{(d): t=8.5e-4}
\includegraphics[width=2in]{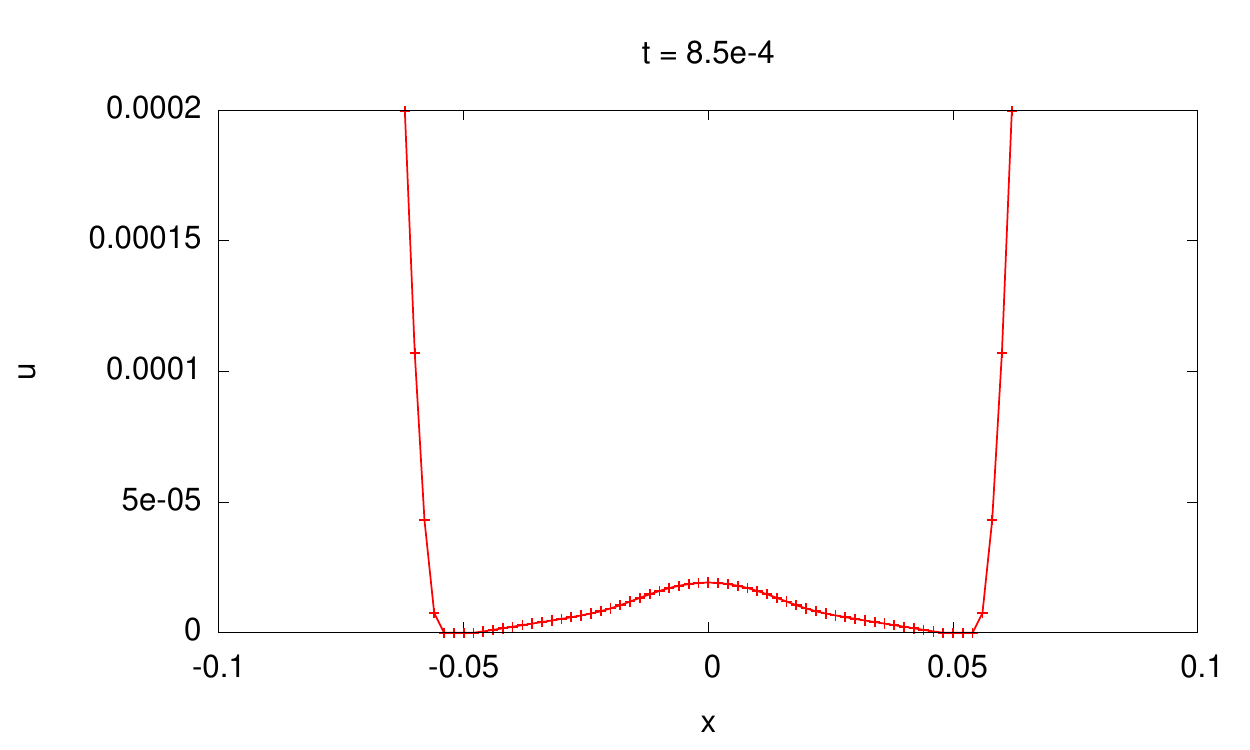}
\end{minipage}
\begin{minipage}[b]{2.0in}
\centerline{(e): t=1.0e-3}
\includegraphics[width=2in]{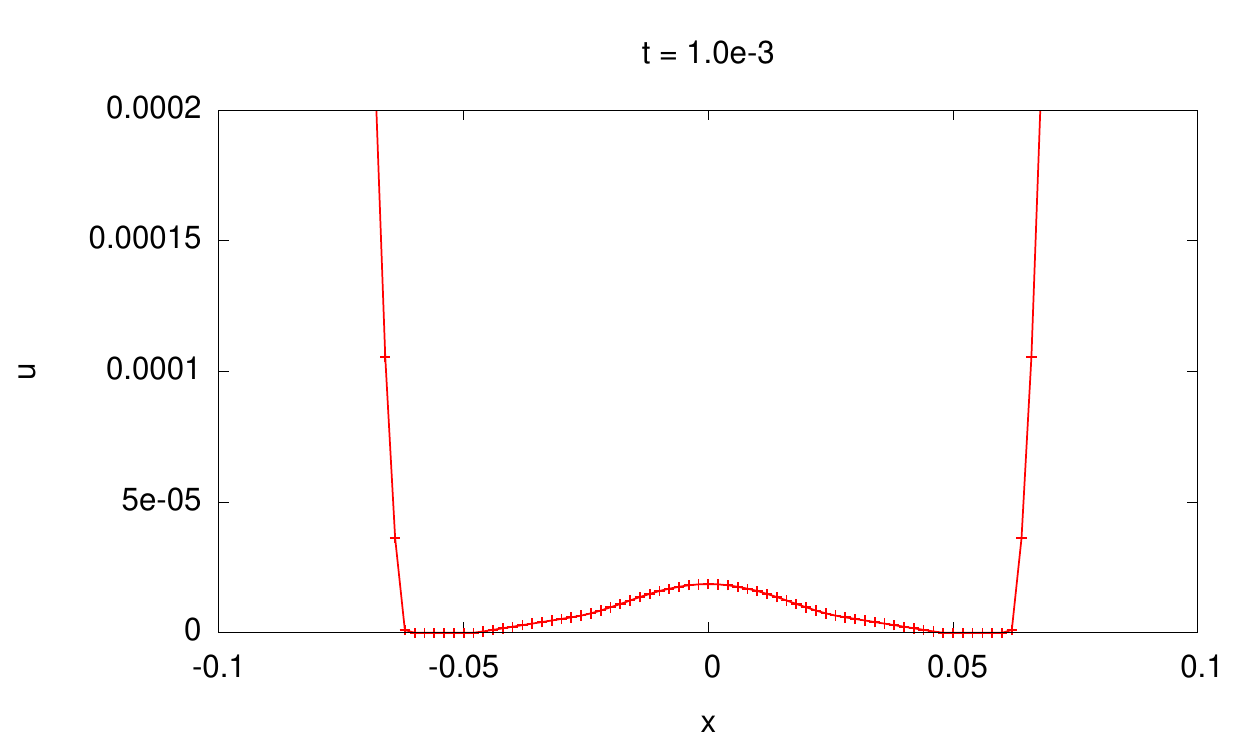}
\end{minipage}
\begin{minipage}[b]{2.0in}
\centerline{(f): t=1.5e-3}
\includegraphics[width=2in]{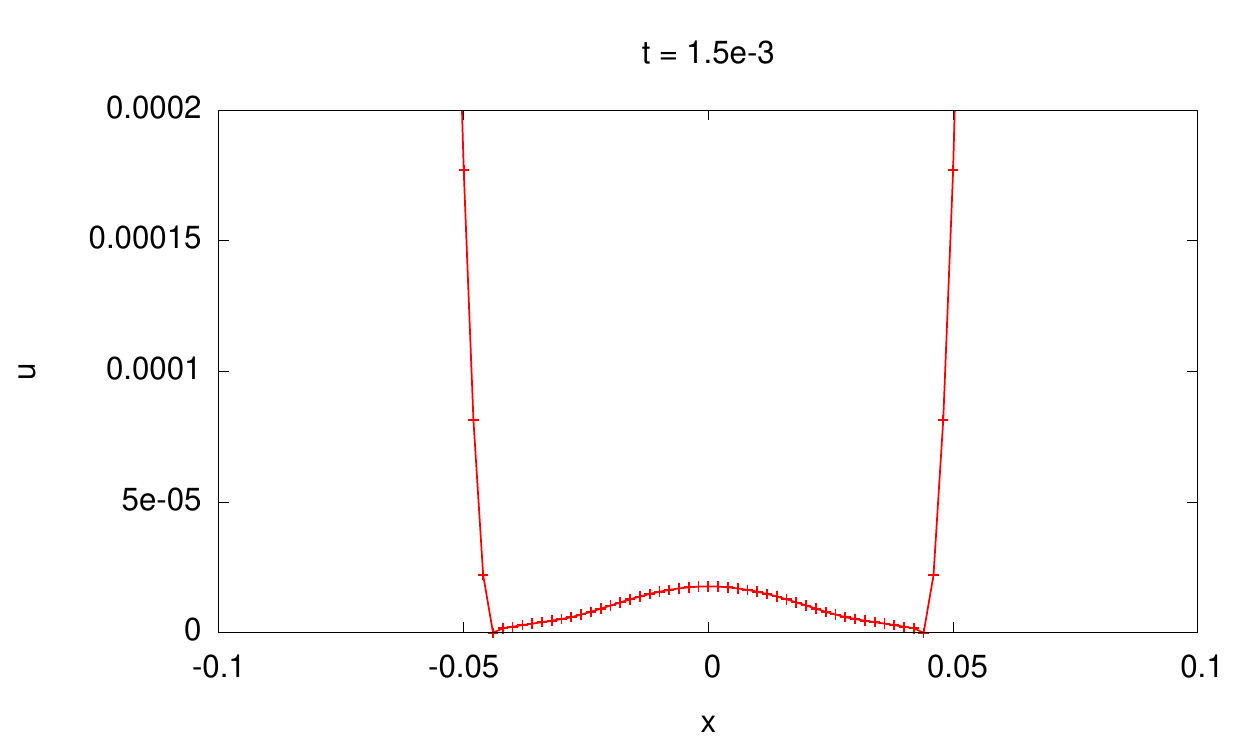}
\end{minipage}
}
\hbox{
\begin{minipage}[b]{2.0in}
\centerline{(g): t=2.0e-3}
\includegraphics[width=2in]{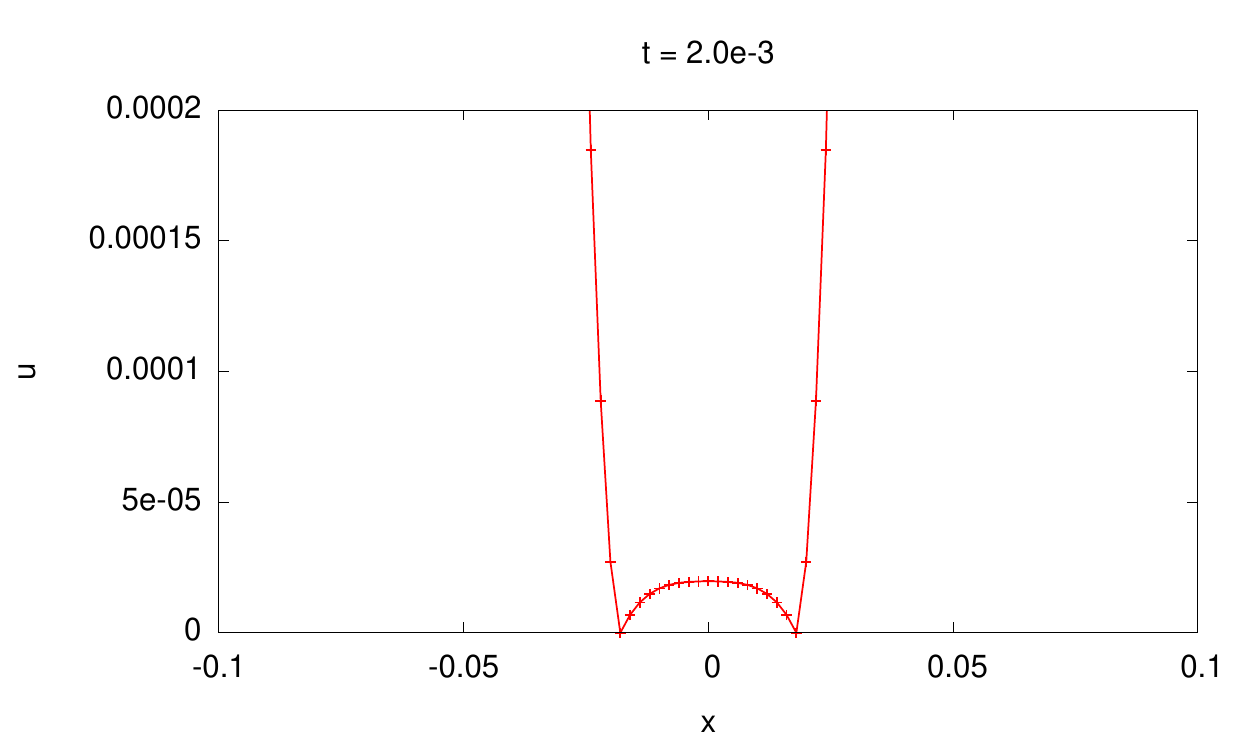}
\end{minipage}
\begin{minipage}[b]{2.0in}
\centerline{(h): t=2.3e-3}
\includegraphics[width=2in]{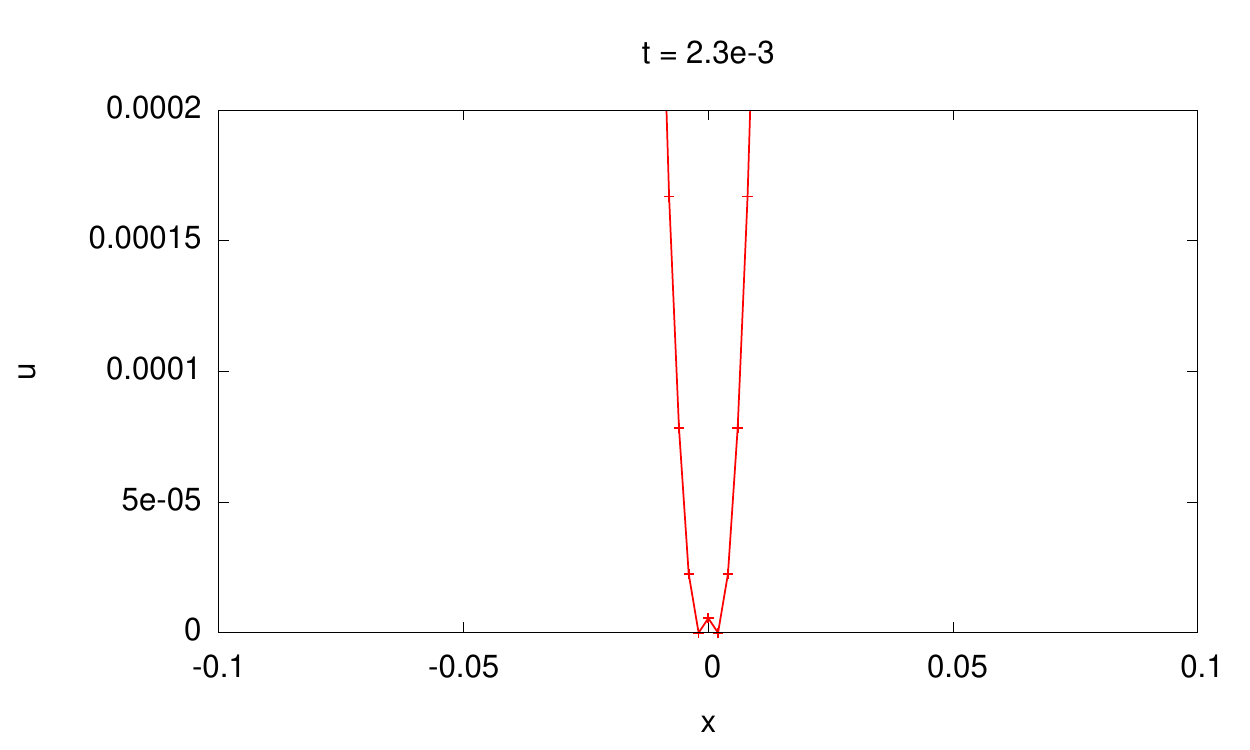}
\end{minipage}
\begin{minipage}[b]{2.0in}
\centerline{(i): t=2.34e-3}
\includegraphics[width=2in]{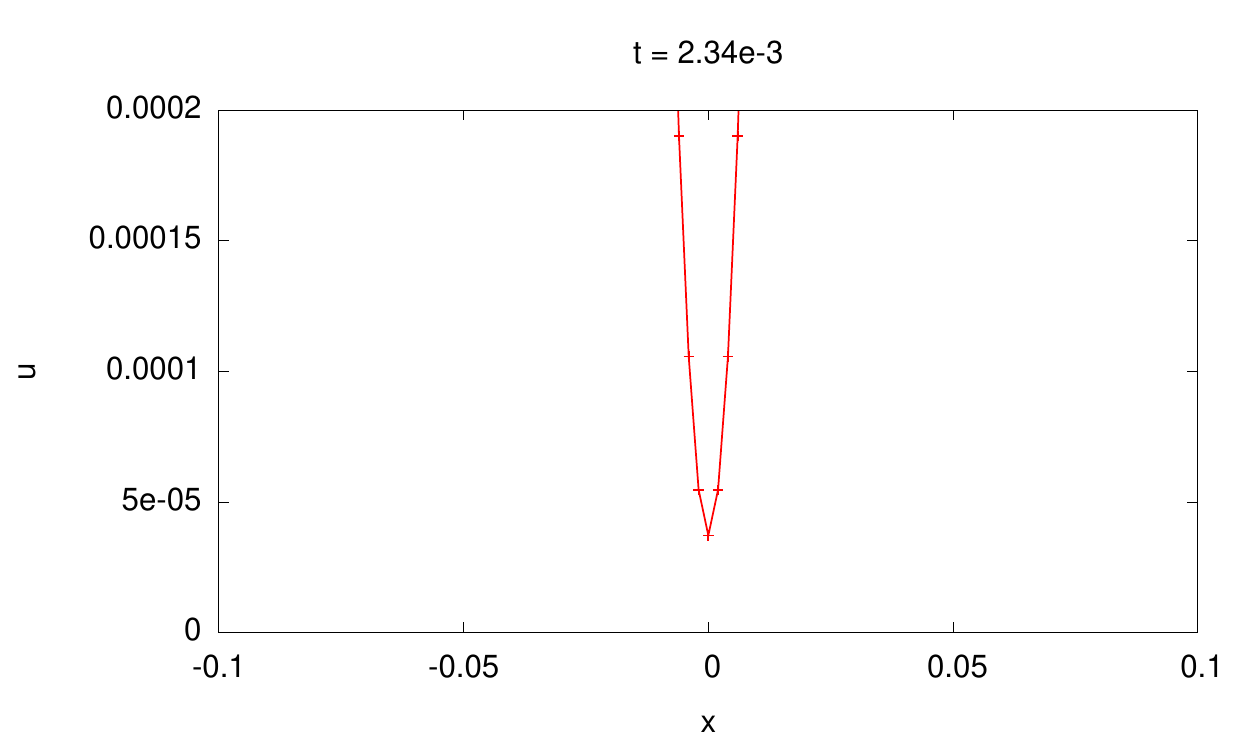}
\end{minipage}
}
\caption{The one-dimensional regularized lubrication-type problem (with $\epsilon = 10^{-14}$).
The close views of the numerical solution
are shown at various time instants during the singularity development. The numerical solution is obtained
by the cutoff method with regularization for $f(u)$ on a uniform mesh of 1001 points with
$\Delta t = 10^{-6}$.}
\label{Exa4.1-4}
\end{figure}

As we have seen in the above, there is no need to use any regularization for $f(u)$ (cf. (\ref{lub-2})) with the cutoff
method. On the other hand, it is interesting to see the effects of the regularization of $f(u)$ on the solution.
We first point out that the regularization does not guarantee the nonnegativity of the solution for
the central finite difference discretization on a uniform mesh. Thus, we also need to use
the cutoff method for the regularized equation (\ref{lub-2}).
A solution obtained for (\ref{lub-2}) (with $\epsilon = 10^{-14}$) with the cutoff method is shown
in Fig.~\ref{Exa4.1-4}. By a direct comparison with Fig.~\ref{Exa4.1-3}, one can see that there
is only a slight difference (at the level of $\mathcal{O}(10^{-6})$) in the onset time and disappearance
time of the singularity. Moreover, the solution of the regularized problem does not become identically
zero on an interval for a period of time. Instead, there is a bump at the central part of the interval and
the height of the bump maintains constant almost for the whole appearance period of singularity.
Once again, the length of the interval first increases from zero, reaches its maximum, and decreases to zero
when the solution becomes positive again. 

\begin{figure}[thb]
\centering
\hbox{
\begin{minipage}[b]{3.0in}
\centerline{(a): Solutions at onset of singularity}
\includegraphics[width=3.0in]{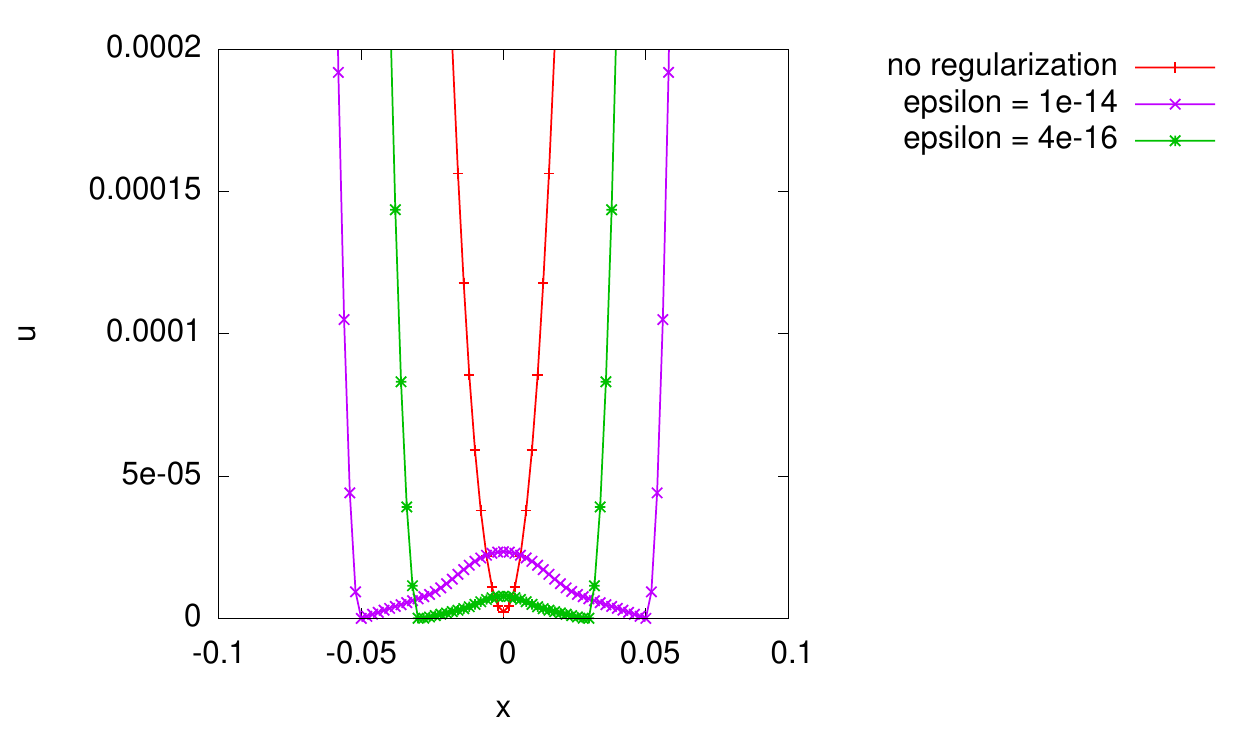}
\end{minipage}
\begin{minipage}[b]{3.0in}
\centerline{(b): Solutions at $t=0.001$}
\includegraphics[width=3.0in]{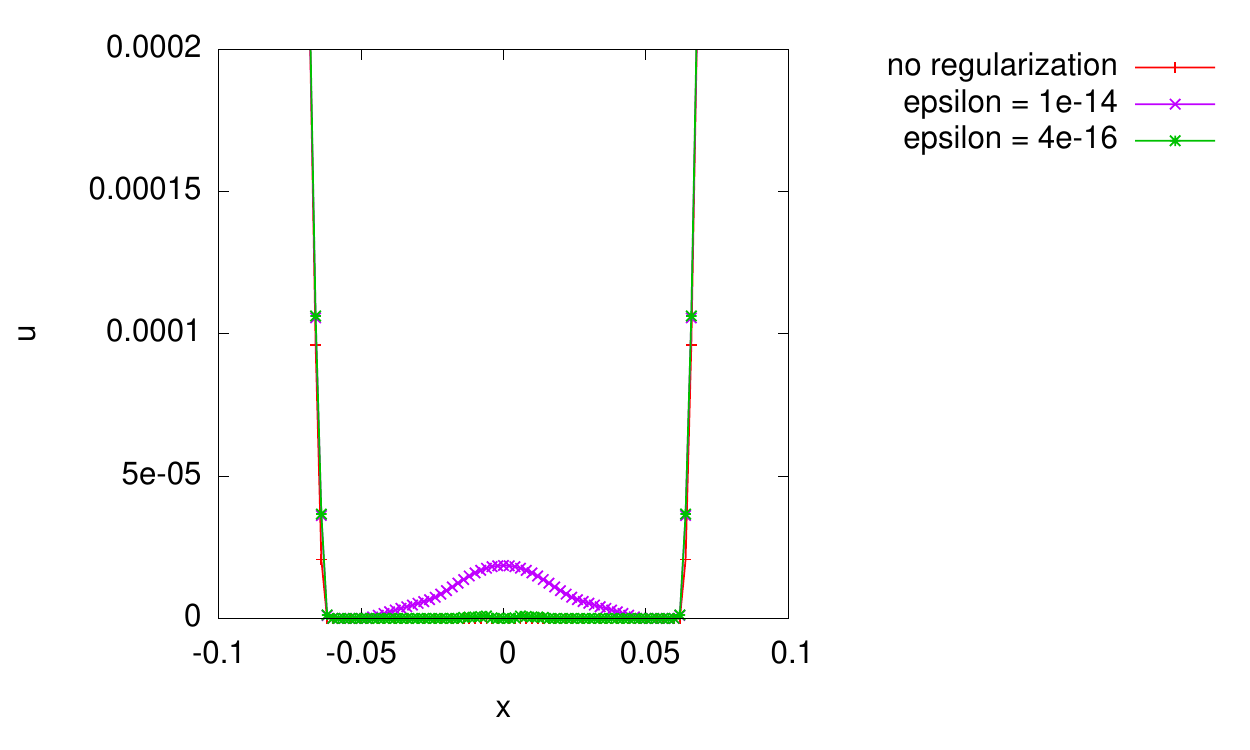}
\end{minipage}
}
\caption{The one-dimensional lubrication-type problem with/without regularization. Numerical solutions
at onset of singularity or at $t=0.001$.}
\label{Exa4.1-6}
\end{figure}

Fig.~\ref{Exa4.1-6} shows comparison of numerical solutions at onset of singularity and at $t= 0.001$.
It can be seen that the regularization changes the onset pattern: the solution to the non-regularized equation
touches the $x$-axis at one point whereas those to the regularized one touch at two points simultaneously.
Moreover, the length of the touching interval and the height of the bump at onset depend on the regularization
parameter $\epsilon$. The smaller $\epsilon$, the smaller the touching interval and the lower the bump.
This suggests $u^\epsilon \to u$ as $\epsilon \to 0$. It is interesting to point out that all solutions, for regularization
or non-regularization, have almost the same length of the touching interval at $t=0.001$.

\begin{figure}[thb]
\centering
\hbox{
\begin{minipage}[b]{3.0in}
\centerline{(a): Numerical solution at $t=0.001$}
\includegraphics[height=2.5in]{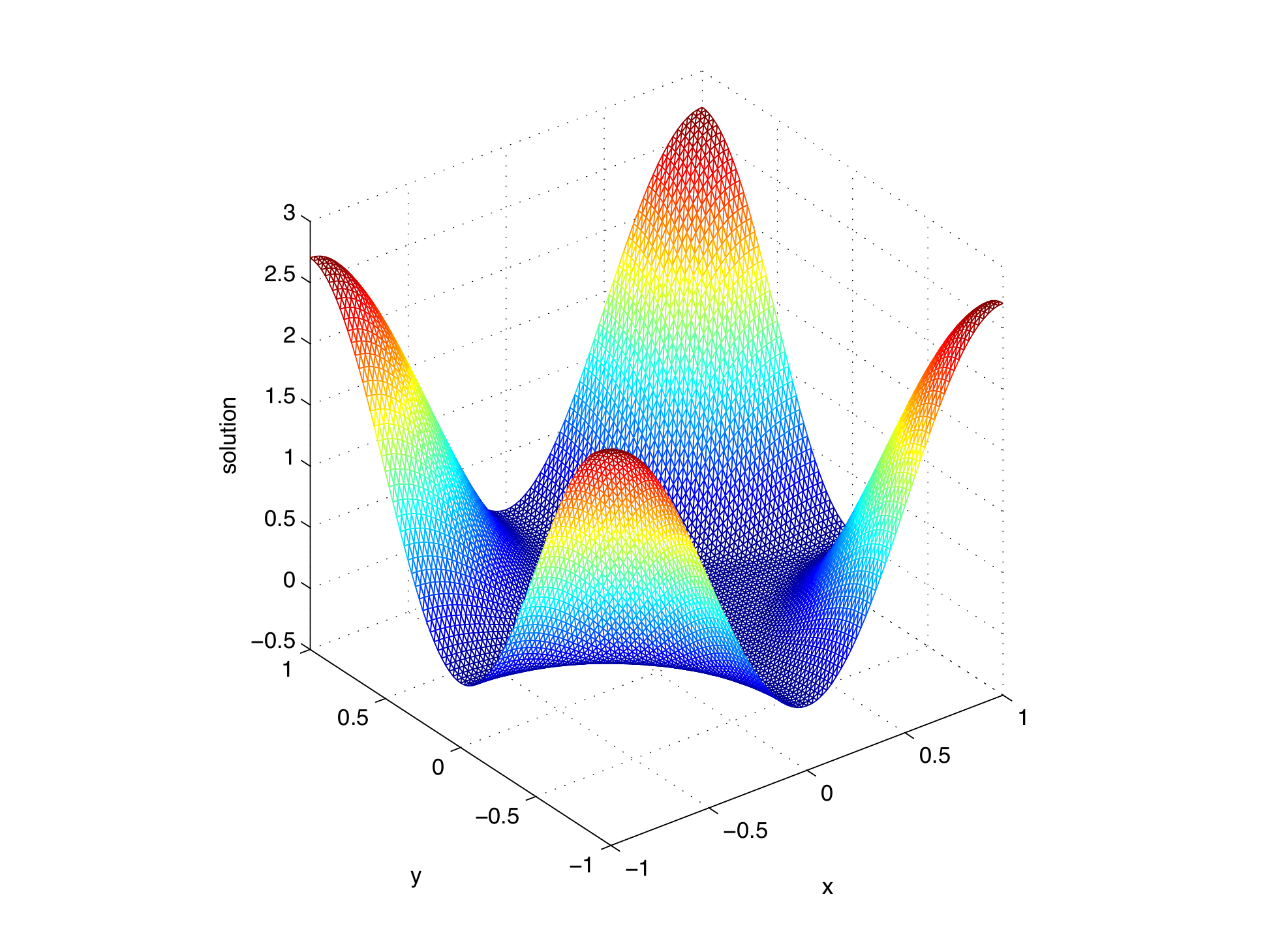}
\end{minipage}
\begin{minipage}[b]{3.0in}
\centerline{(b): Solution along $x=y$ at various times}
\includegraphics[height=2.5in]{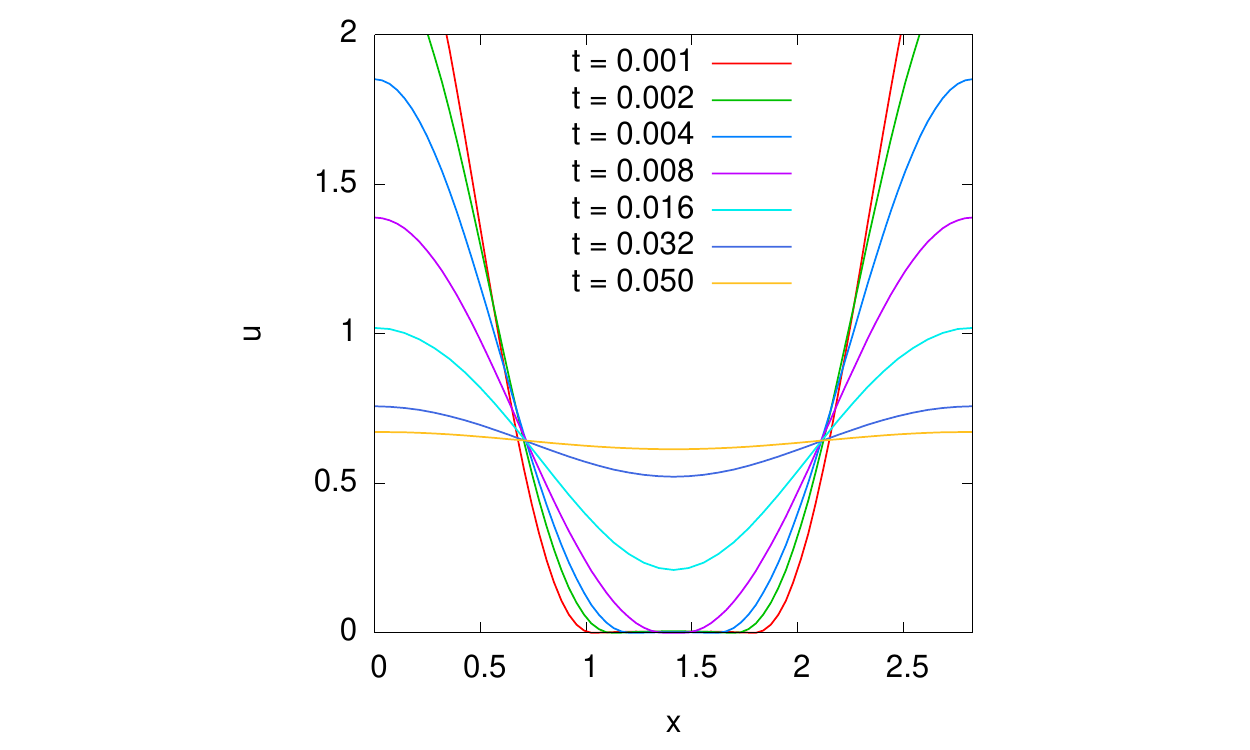}
\end{minipage}
}
\caption{The two-dimensional lubrication problem. The solution is obtained with the cutoff method with
a uniform mesh of size $81\times 81$.}
\label{Exa4.2-1}
\end{figure}

\begin{figure}[thb]
\centering
\hbox{
\begin{minipage}[b]{3.0in}
\centerline{(a): $U^N$}
\includegraphics[width=3.5in]{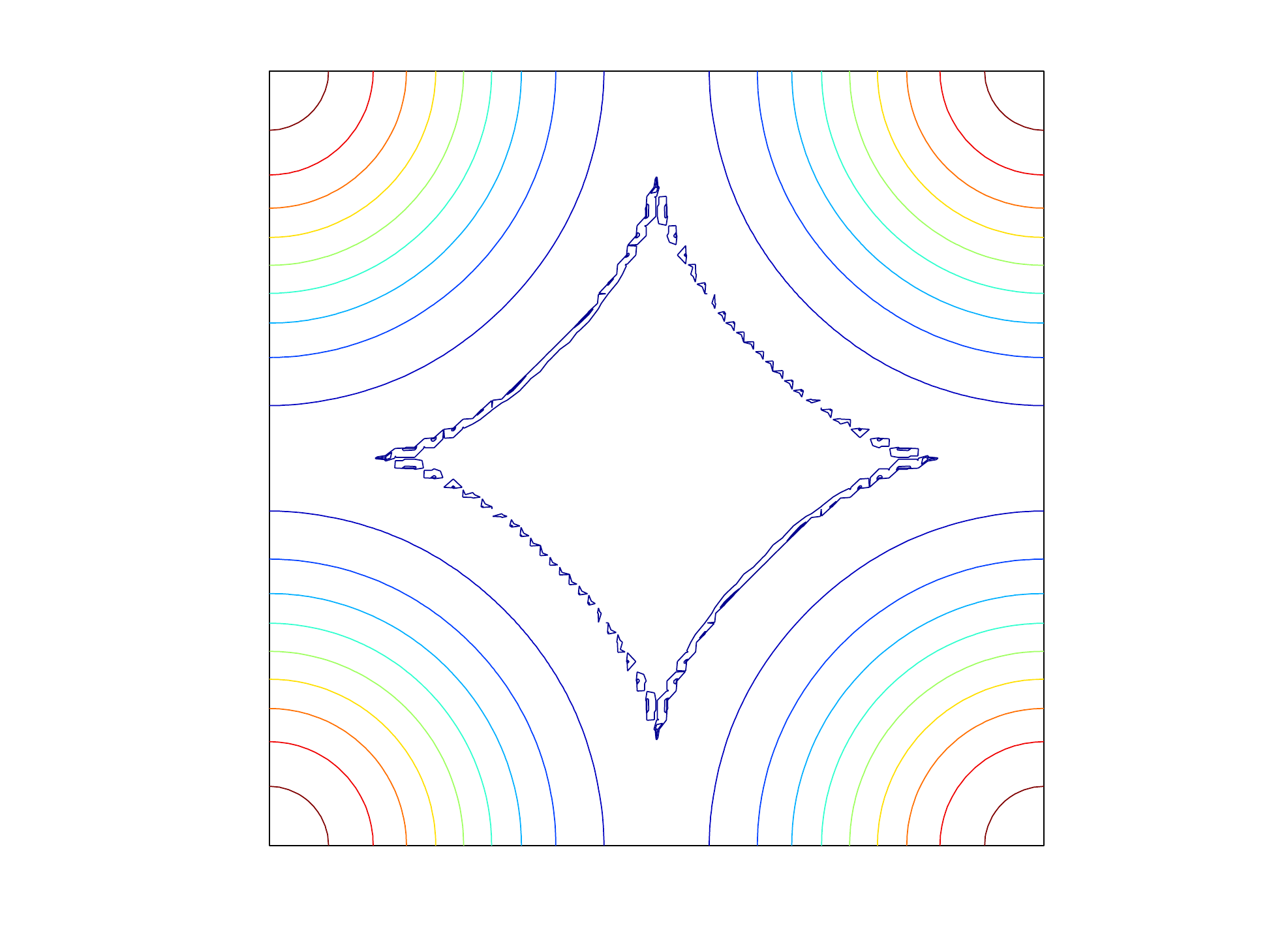}
\end{minipage}
\begin{minipage}[b]{3.0in}
\centerline{(b): $(U^N)^+$}
\includegraphics[width=3.5in]{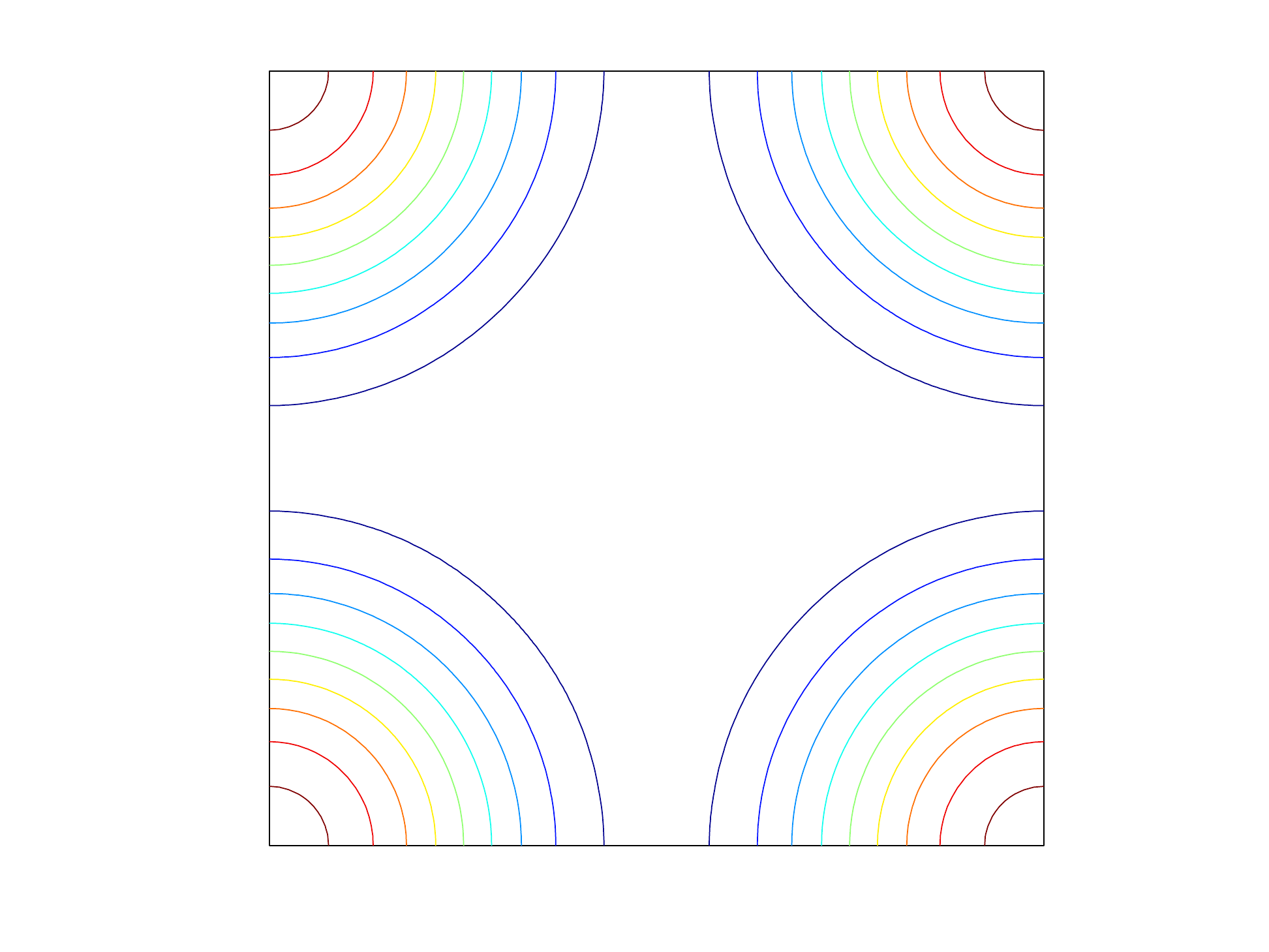}
\end{minipage}
}
\caption{The two-dimensional lubrication problem. Contours of the numerical solutions at $t=0.001$ before and after cutoff.}
\label{Exa4.2-2}
\end{figure}

Finally, we present some numerical results for an IBVP of the two-dimensional lubrication-type equation as
\beq
\begin{cases}
& u_t + \nabla \cdot (f(u) \nabla \Delta u) = 0, \quad  (x,y) \in  \Omega \equiv [-1,1] \times [-1,1], \;\; f(u) = u^{\frac{1}{2}},
\\
& u = 0, \; \; \frac{\p \Delta u}{\p n}  = 0, \quad (x,y) \in  \partial \Omega
\\
& u(0, x,y) =( 0.8 -\cos (\pi x) + 0.25 \cos (2\pi x) ) ( 0.8 -\cos (\pi y) + 0.25 \cos (2\pi y) ), \quad  (x,y) \in  \Omega .
\end{cases}
\label{ibvp-5}
\eeq
Fig.~\ref{Exa4.2-1} (a) shows a numerical solution at $t = 0.001$ obtained with a uniform mesh of size $81 \times 81$
while its profiles along the diagonal line $x=y$ at various times are shown in Fig.~\ref{Exa4.2-1} (b).
One can see that the development of singularity has a similar behavior as in one dimension (cf. Fig.~\ref{Exa4.1-1}).
One may also see that the onset time for the current example is slightly earlier. This is because the current
initial solution has a smaller minimum value 0.025.
The contours of the numerical solutions at $t=0.001$ before and after cutoff are shown
in Fig.~\ref{Exa4.2-2}. The undershoot (with $u_{min}=-3.23\times 10^{-4}$) is visible
in the solution before cutoff. These results demonstrate that the cutoff method can be used
for the simulation of singularity development in the two-dimensional lubrication-type problem
without need of any type of regularization.

\section{Conclusions and further comments}
\label{SEC:conclusion}

In the previous section we have studied the cutoff method for the numerical computation of
nonnegative solutions of parabolic partial differential equations. Several properties
of the cutoff method are given in Lemmas~\ref{lem2.1} -- \ref{lem2.3}. Convergence of a class of
finite difference methods is proved (Theorems \ref{thm2.1} and \ref{thm2.2}) when they are incorporated
with the cutoff method for linear parabolic PDEs. The method is investigated for two applications, linear anisotropic
diffusion problems and nonlinear lubrication-type PDEs. The numerical results are consistent with theoretical
predictions and in good agreement with existing results in the literature.

\begin{figure}[thb]
\centering
\hbox{
\begin{minipage}[b]{3in}
\centerline{(a): $U^N$}
\includegraphics[width=3.5in]{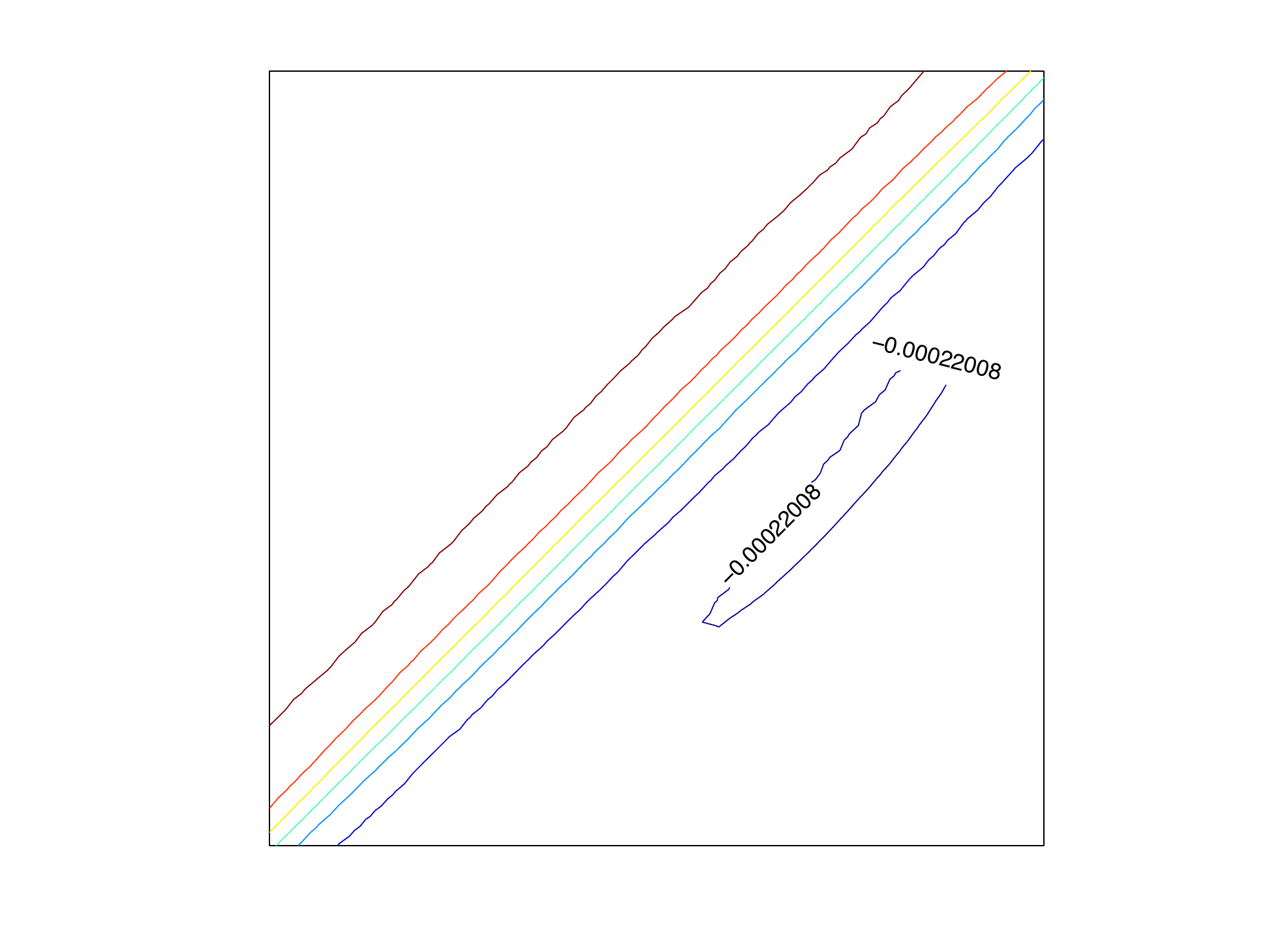}
\end{minipage}
\begin{minipage}[b]{3in}
\centerline{(b): $(U^N)^+$}
\includegraphics[width=3.5in]{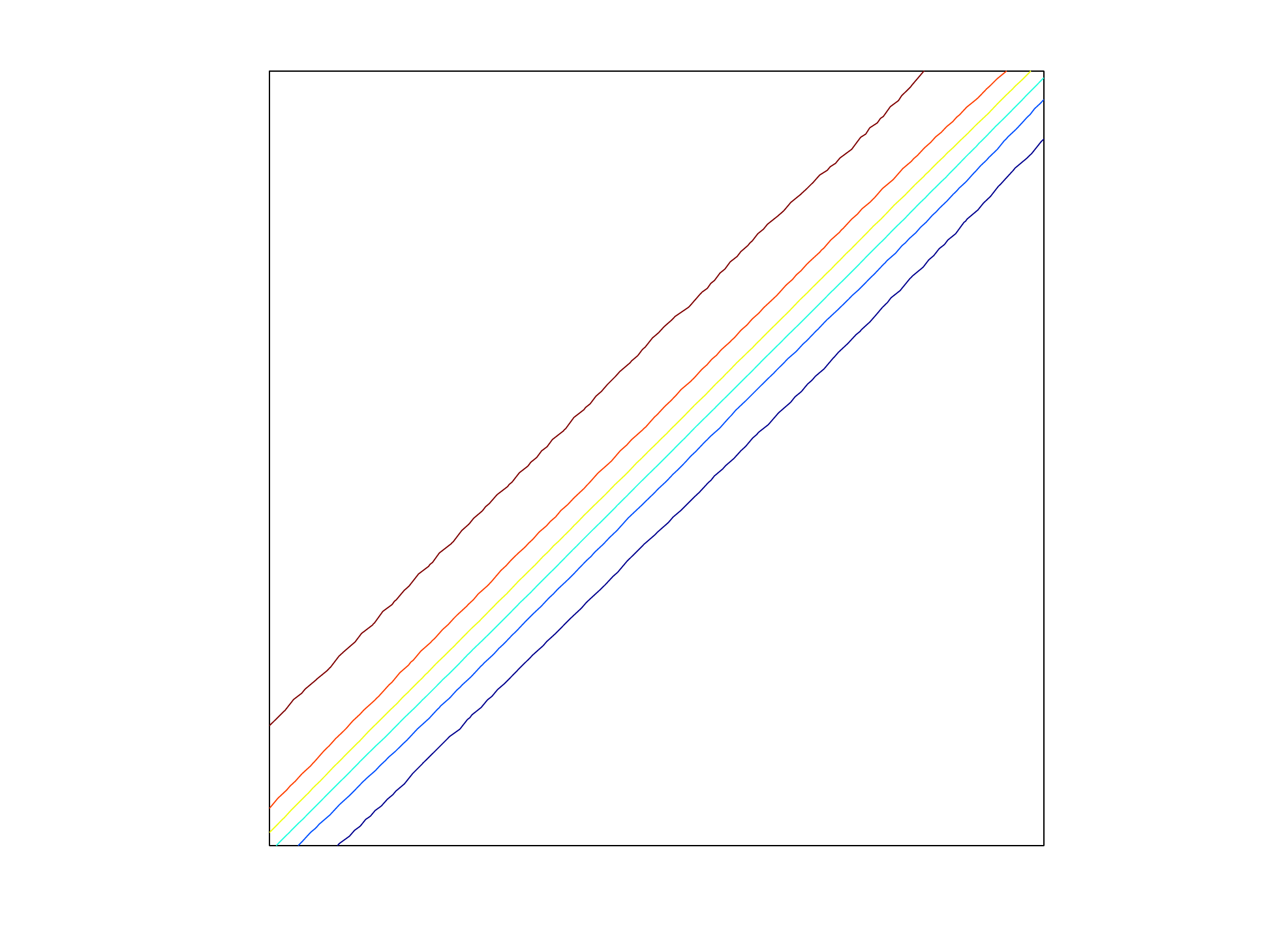}
\end{minipage}
}
\caption{The standard linear finite element method for Example~(\ref{ibvp-2}).
Contours of the numerical solutions at $t = 1$ before and after cutoff obtained 
with a Delaunay mesh of $N=2350$ elements.}
\label{Exa3.1-6}
\end{figure}

\begin{figure}[thb]
\centering
\includegraphics[scale=0.8]{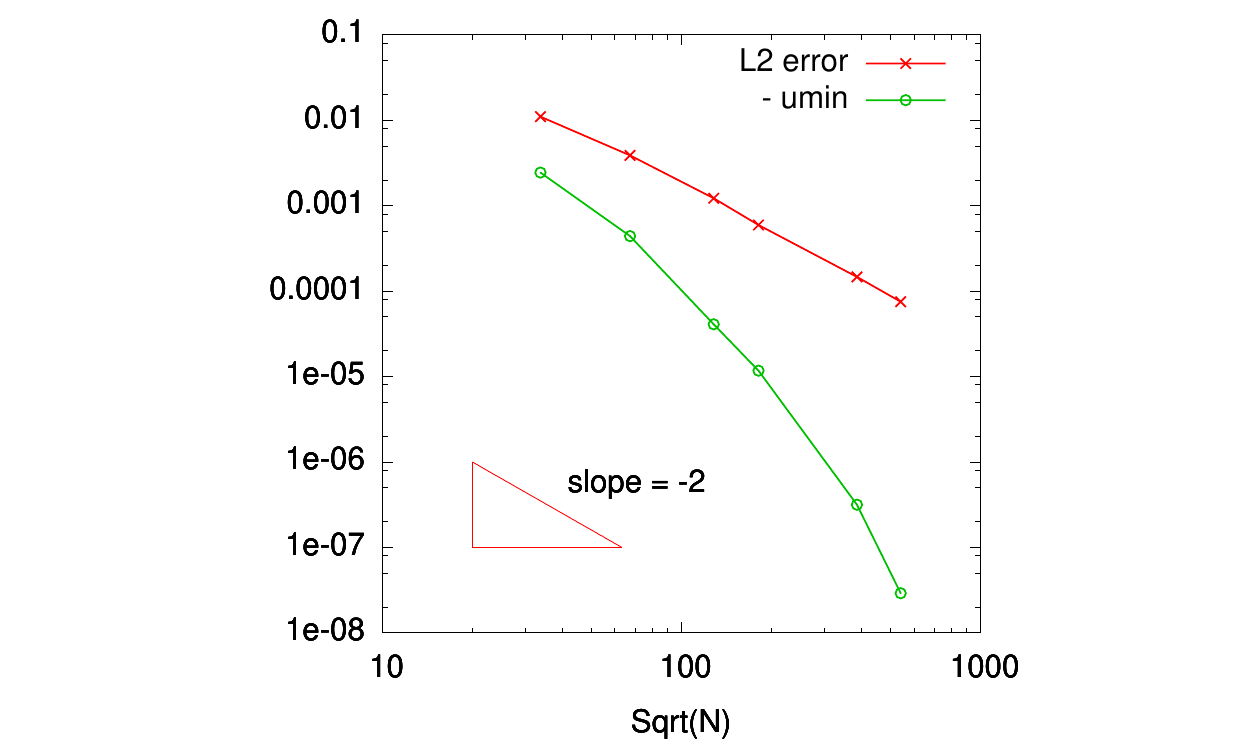}
\caption{The standard linear finite element method for Example~(\ref{ibvp-2}).
The $L^2$ error, $\| (U^N)^{+}-u^N \|_{L^2(\Omega)}$, and the maximal undershoot, $-u_{min}$,
are shown as functions of $\sqrt{N}$ with $N$ being the number of elements.
}
\label{Exa3.1-7}
\end{figure}

We have considered finite difference methods in this work.
But it is worth pointing out that the cutoff method can also be used with other discretization methods, e.g.
collocation, finite volume, finite element, or spectral methods.
As an example, we show in Figs.~\ref{Exa3.1-6} and \ref{Exa3.1-7} results obtained with the standard linear
finite element method on Delaunay meshes for anisotropic diffusion problem (\ref{ibvp-2}).
They are comparable with those in Figs.~\ref{Exa3.1-2} and \ref{Exa3.1-3}.
Theoretical analysis for finite element methods with the cutoff strategy is currently underway.


\vspace{20pt}

{\bf Acknowledgment.} 
The work was supported in part by the National
Science Foundation (U.S.A.) under Grants DMS-1115118 and DMS-1115408
and the Natural Science Foundation of China under Grant 40906048.


\end{document}